\documentclass[11pt]{article} 
\usepackage{a4,amssymb, amsmath, amsthm}

\usepackage{a4,amssymb, amsmath}
\usepackage{graphics}

\newtheorem{theorem}{Theorem}[section]

\newtheorem{lemma}[theorem]{Lemma}

\newtheorem{definition}[theorem]{Definition}
\newtheorem{remark}[theorem]{Remark}

\begin{document}

\title{A priori error estimates for a time-dependent boundary element method for the acoustic wave equation in a half-space}
\author{Heiko Gimperlein, Zouhair Nezhi and Ernst P.~Stephan}

\date{}

\maketitle

\begin{abstract}
We investigate a time--domain Galerkin boundary element method for the wave equation outside a Lipschitz obstacle in an absorbing half--space, with application to the sound radiation of tyres. A priori estimates are presented for both closed surfaces and screens, and we discuss the relevant properties of anisotropic Sobolev spaces and the boundary integral operators between them.
\end{abstract}

\vskip 1.0cm

\section{Introduction}\label{intro}

Motivated by the sound radiation of tyres {on} a street, this article provides the analytical background to analyze a time--domain boundary element method for the direct scattering problem for the wave equation outside an obstacle in an absorbing half--space.\\

{Let $d \geq 2$ and} ${\Omega^i}\subset \mathbb{R}^d_+$ be a bounded Lipschitz domain such that the exterior domain $\Omega^e =\mathbb{R}^d_+ \backslash \overline{\Omega^i}$ is Lipschitz and connected.  {The reader may wish to think of $\Omega^i$ as a solid tyre, either in contact with the street (on $\partial \Omega^i \cap \partial \mathbb{R}^d_+$) or elevated above it ($\partial \Omega^i \cap \partial \mathbb{R}^d_+ = \emptyset$ ).} The boundary of $\Omega^e$ decomposes into the boundary $\Gamma = \partial \Omega^e \cap \partial \Omega^i$ of the obstacle and the boundary $\Gamma_\infty = \partial \Omega^e \cap \partial\mathbb{R}^d_+$ of the half--space. In general, $\Gamma$ is a Lipschitz manifold with boundary, and we {emphasize} the case $d=3$.\\

We aim to find a weak solution to an acoustic initial boundary problem for the wave equation in $\Omega^e$:
\begin{align}\label{eq:strong_half_space}
\frac{\partial^2 u}{\partial t^2} -\Delta u&=0 \quad\mbox{ in } \mathbb{R}^+\times \Omega^e \nonumber\\
u(0,x)=\frac{\partial u}{\partial t}(0,x)&=0   \quad\mbox{ in } \Omega^e\\
\frac{\partial u}{\partial n} -\alpha \frac{\partial u}{\partial t} &= g\quad\mbox{ on }\mathbb{R}^+\times \Gamma  \nonumber\\
\frac{\partial u}{\partial n} -\alpha_\infty \frac{\partial u}{\partial t} &= 0\quad\mbox{ on } \mathbb{R}^+\times \Gamma_\infty \nonumber\,.
\end{align}
Here $n$ denotes the inward unit normal vector to $\partial \Omega^e$, $g$ lies in a suitable Sobolev space, $\alpha\in L^\infty(\Gamma)$ and $\alpha_\infty \in \mathbb{C}$. We also consider the simpler Dirichlet problem on $\Gamma$, for which instead {of the absorbing boundary condition,} $u|_{\mathbb{R}^+\times\Gamma}$ is given.\\

This article reduces the acoustic and Dirichlet boundary problems to time--dependent integral equations on $\mathbb{R}^+\times \Gamma$ and studies a Galerkin time--domain boundary element method for their approximation. {Time--dependent Galerkin boundary element methods for wave problems were introduced by Bamberger and Ha-Duong \cite{bamberger86}. Some relevant works on the numerical implementation of the resulting marching-in-on-time scheme include the Ph.D.~thesis of Terrasse \cite{terrasse93} and \cite{hld03}, with fast methods developed in the engineering literature \cite{mic}. Alternative ansatz functions in time have been explored in \cite{dd1, dd2}. A detailed exposition of the mathematical background of time--domain integral equations and their discretizations is available in the lecture notes by Sayas \cite{sayas}. \\
In the special case of the half-space, our work is motivated by the recent explicit formulas for the fundamental solutions obtained by Ochmann \cite{Ochmann01}, which include acoustic boundary conditions on the surface of the street.}\\

Section 2 introduces space--time anisotropic Sobolev spaces of supported resp.~extendable distributions on $\mathbb{R}^+\times \Gamma$. Their approximation theory and interpolation operators are the subject of Section 3. Subsequent sections follow the approach of Bamberger and Ha Duong \cite{bamberger86, Ha-Duong03a}, see also \cite{costabel04}, to analyze the coercivity and boundedness properties of time--dependent layer potentials adapted to the acoustic boundary conditions on $\Gamma_\infty$. As conclusion, we deduce a priori error estimates for the Galerkin solutions.\\

The results in this article provide a basic theoretical background for further theoretical and computational analysis. Based on the set--up presented here, future work will address a posteriori error estimates and adaptive procedures {\cite{Gimperlein}} as well as numerical studies of engineering benchmarks {\cite{Banz}}.\\

We acknowledge support within the project ``LeiStra3'' by the German Bundesministerium f\"{u}r Wirtschaft und Energie as well as the Bundesanstalt f\"{u}r Stra\ss enwesen. H.G.~thanks the Danish Science Foundation (FNU) for partial support through research grant 10-082866.\\

{\emph{Notation:} To simplify notation, we will write $f \lesssim g$, if there exists a constant $C>0$ independent of the arguments of the functions $f$ and $g$ such that $f \leq C g$. We will write $f \lesssim_\sigma g$, if $C$ may depend on $\sigma$.}

\section{Space--time anisotropic Sobolev spaces}

Space--time anisotropic Sobolev spaces on the boundary $\Gamma$ provide a convenient setting to study the mapping properties of the layer potentials \cite{Ha-Duong03a, costabel04}.
We define {these Sobolev spaces} on a general $C^{k-1,1}$ closed, orientable manifold $\mathcal{M}$ with boundary, in particular {obtain Sobolev spaces} on $\mathcal{M}=\Gamma$ and on $\mathcal{M}=\Omega$. {The isotropic case is well-known from elliptic problems, see \cite{tri}, and \cite{screen} when $\partial\mathcal{M} \neq \emptyset$.} \\

If $\partial\mathcal{M}\neq \emptyset$, first extend $\mathcal{M}$ to a $C^{k-1,1}$, closed, orientable manifold $\widetilde{\mathcal{M}}$. For example, if $\mathcal{M} = \partial\Omega \cap \mathbb{R}^d_+$ and $ k=1$, we may take $\widetilde{\mathcal{M}}$ to be union of $\partial\Omega \cap\overline{\mathbb{R}^d_+}$ and its image under reflection at $\partial \mathbb{R}^d_+$.

On $\mathcal{M}$, the usual Sobolev spaces of supported distributions may now be considered {for $r \in \mathbb{R}$}:
$$\widetilde{H}^r(\mathcal{M}) = \{u\in H^r(\widetilde{\mathcal{M}}): \mathrm{supp}\ u \subset {\overline{\mathcal{M}}}\}\ .$$
${H}^r(\mathcal{M})$ is defined as the quotient space $ H^r(\widetilde{\mathcal{M}}) / \widetilde{H}^r({\widetilde{\mathcal{M}}\setminus\overline{\mathcal{M}}})$. \\
To define a family of Sobolev norms, introduce a partition of unity subordinate to the covering of $\widetilde{\mathcal{M}}$ by open sets $B_i$. For a partition of unity $\alpha_i$ and diffeomorphisms $\varphi_i$ mapping each $B_i$ into the unit cube $Q\subset \mathbb{R}^n$, a family of equivalent Sobolev norms is induced from $\mathbb{R}^n$:
\begin{equation*}
 ||u||_{r,\omega,{\widetilde{\mathcal{M}}}}=\left( \sum_{i=1}^p \int_{\mathbb{R}^n} (|\omega|^2+|\xi|^2)^r|\mathcal{F}\left\{(\alpha_i u)\circ \varphi_i^{-1}\right\}(\xi)|^2 d\xi \right)^{\frac{1}{2}}\ .
\end{equation*}
Here $\omega \in \mathbb{C}\setminus \{0\}$ and $\mathcal{F}$ denotes the Fourier transform. {This norm on $H^r(\widetilde{\mathcal{M}})$ induces a norm on $H^r(\mathcal{M})$ as $||u||_{r,\omega,\mathcal{M}} = \inf_{v \in \widetilde{H}^r(\widetilde{\mathcal{M}}\setminus\overline{\mathcal{M}})} \ ||u+v||_{r,\omega,\widetilde{\mathcal{M}}}$.} \\

{The weighted norm on $\widetilde{H}^r(\mathcal{M})$ is defined as $||u||_{r,\omega,\mathcal{M}, \ast } = ||e_+ u||_{r,\omega,\widetilde{\mathcal{M}}}$, where $e_+$ extends the distribution $u$ by $0$ from $\mathcal{M}$ to $\widetilde{\mathcal{M}}$. It is stronger than $||u||_{r,\omega,\mathcal{M}}$ whenever $r \in \frac{1}{2} + \mathbb{Z}$.}\\

For $|r|\leq k$ the thus defined Sobolev spaces are independent of the choice of $\alpha_i$ and $\varphi_i$. \\

Using these norms, the trace operator from $H^1(\Omega)$ to $H^{\frac{1}{2}}(\Gamma)$ is continuous in the $\omega$--dependent norms. As shown in \cite{Ha-Duong03a}, for $\sigma >0$ its operator norm is uniformly bounded in the half-plane $\{\omega \in \mathbb{C} : \mathrm{Im} \ \omega > \sigma\}$ by a function of $\sigma$ alone. It admits a right inverse from $H^{\frac{1}{2}}(\Gamma)$ and $\widetilde{H}^{\frac{1}{2}}(\Gamma)$ to $H^1(\Omega)$ whose norm is similarly bounded in terms of $\sigma$.\\

Let $E$ be a Hilbert space. We define\\
\begin{equation*}
 LT(\sigma,E) = \{ f \in \mathcal{D}^{'}_{+}(E); e^{-\sigma t} f \in \mathcal{S}^{'}_{+}(E)\}\ ,
\end{equation*}
where $\mathcal{D}^{'}_{+}(E)$ resp.~$\mathcal{S}^{'}_{+}(E)$ denote the sets of distributions resp.~tempered distributions on $\mathbb{R}$ with values in $E$ and support in $[0,\infty)$. As $LT(\sigma,E)\subset LT(\sigma',E)$ if $\sigma < \sigma'$, we may define $\sigma(f)=\inf\{\sigma : f \in LT(\sigma,E)\}$.\\
The set of Laplace transformable distributions with values in $E$ is denoted by\\
\begin{equation*}
 LT(E) = \underset{\sigma \in \mathbb{R}}{\cup} LT(\sigma,E)\ .
\end{equation*}
For  $f \in LT(E)$, its Fourier-Laplace transform $\hat{f}(\omega)=\mathcal{F} f(\omega)$ is defined in the complex half--plane $\{\omega\in\mathbb{C} : \mathrm{Im}\ \omega>\sigma(f)\}$.

We recall the well-known Parseval identity in this setting:
\begin{lemma}
For  $f, g \in L^1_{loc}(\mathbb{R},E) \cap LT(E)$ and $\sigma > max(\sigma(f),\sigma(g))$ there holds
\begin{equation*}
 \frac{1}{2 \pi} \int_{\mathbb{R} +i \sigma} (\hat{f}(\omega),\hat{g}(\omega))_E d\omega = \int_{\mathbb{R}} e^{-2\sigma t} (f(t),g(t))_E dt\ .
\end{equation*}
\end{lemma}
Space--time anisotropic Sobolev spaces on $\mathcal{M}$ are now defined as follows:
\begin{definition}\label{sobdef}
For $s, r \in\mathbb{R}$ define
\begin{align*}
 H^s_\sigma(\mathbb{R}^+,{H}^r(\mathcal{M}))&=\{ u \in LT({H^r}(\mathcal{M}))  :   ||u||_{s,r,\mathcal{M}} < \infty \}\, \\
 H^s_\sigma(\mathbb{R}^+,\widetilde{H}^r(\mathcal{M}))=&\{ u \in LT(\widetilde{H}^r(\mathcal{M}))  :   ||u||_{s,r,\mathcal{M},\ast} < \infty \}\,,
\end{align*}
where
\begin{align*}
\|u\|_{s,r,\mathcal{M}}&=\left(\int_{-\infty+i\sigma}^{+\infty+i\sigma}|\omega|^{2s}\ \|\hat{u}(\omega)\|^2_{r,\omega,\mathcal{M}}\ d\omega \right)^{\frac{1}{2}}\ ,\\
\|u\|_{s,r,\mathcal{M},\ast}&=\left(\int_{-\infty+i\sigma}^{+\infty+i\sigma}|\omega|^{2s}\ \|\hat{u}(\omega)\|^2_{r,\omega,\mathcal{M},\ast}\ d\omega \right)^{\frac{1}{2}}\,.
\end{align*}
\end{definition}
As above, the spaces are invariantly defined whenever $|r|\leq k$.

\section{Discretisation}

{For simplicity of notation, in this section we restrict ourselves to the two-- and three--dimensional cases, $d=2$ or $3$.} If $\Gamma$ is not polygonal we approximate it by a piecewise polygonal {curve resp.~}surface and write $\Gamma$ again for the approximation. For simplicity, {when $d=3$} we will use here a surface composed of $N$ triangular facets $\Gamma_i$ such that $\Gamma=\cup_{i=1}^N \Gamma_i$. {When $d=2$, we assume $\Gamma=\cup_{i=1}^N \Gamma_i$ is composed of line segments $\Gamma_i$.} In each case, the elements $\Gamma_i$ are closed  with $int(\Gamma_i) \neq \varnothing$, and for distinct $\Gamma_i,\ \Gamma_j \subset \Gamma$ the intersection $int(\Gamma_i) \cap int(\Gamma_j) = \varnothing$. \\

For the time discretisation we consider a uniform decomposition of the time interval $[0,\infty)$ into subintervals $I_n=[t_{n-1}, t_n)$ with time step $|I_n|=\Delta t$, such that $t_n=n\Delta t \; (n=0,1,\dots)$.\\
We choose a basis ${ \varphi_1^p , \cdots , \varphi_{N_s}^p }$ of the space $V_h^p$ of piecewise polynomial functions of degree $p$ in
space (continuous and vanishing at $\partial \Gamma$ if $p\geq 1$) and  a basis ${\beta^{1,q},\cdots,\beta^{N_t,q}}$  of  the space $V^q_{\Delta t}$ of piecewise  polynomial  functions of degree of $q$ in time (continuous and vanishing at $t=0$ if $q\geq 1$).\\
Let $\mathcal{T}_S=\{T_1,\cdots,T_{N_s}\}$ be the spatial mesh for $\Gamma$ and $\mathcal{T}_T={[0,t_1),[t_1,t_2),\cdots,[t_{N_t-1},T)}$ the time mesh for a finite subinterval $[0,T)$.\\
We consider the tensor product of the approximation spaces in space and time, $V_h^p$ and $V^q_{\Delta t}$, associated to the space--time mesh $\mathcal{T}_{S,T}=\mathcal{T}_S \times\mathcal{T}_T$, and we write
\begin{equation*}
 V^{p,q}_{h,\Delta t}= V_h^p \otimes V^q_{\Delta t}\,.
\end{equation*}
In this section we discuss the projection operators onto $ V^{p,q}_{h,\Delta t}$ and their approximation properties.
We recall the well-known results for $V_h^p$ and $V^q_{\Delta t}$, which we are going to need:
\begin{lemma}\label{proj_time}
Let $\varPi_{\Delta t}$ the orthogonal projection from $L^2(\mathbb{R}_+)$ to $V^q_{\Delta t}$ and $m \leq q$. Then for $s \in [-\frac{1}{2}, \frac{1}{2}]$
\begin{equation*}
 ||f-\varPi_{\Delta t}f||_{\sigma,s,\mathbb{R}_+}\leq C_k \Delta t^{q+1-s} |f|_{\sigma,q+1, \mathbb{R}_+}\,.
\end{equation*}
\end{lemma}
\begin{lemma}\label{proj_space}
Let $\varPi_h$ the orthogonal projection from $L^2(\Gamma)$ to $V_h^p$ and $m\leq p$. Then for $s \in [-\frac{1}{2}, \frac{1}{2}]$ we have in the norms of $H^{s}(\Gamma)$ resp.~$\widetilde{H}^{s}(\Gamma)$:
\begin{align*}
 ||f-\varPi_hf||_{s, \Gamma} &\leq C h^{m+1-s} |f|_{m+1, \Gamma}\\
 ||f-\varPi_hf||_{s, \Gamma,\ast} &\leq C h^{m+1-s} |f|_{m+1, \Gamma}
\end{align*}
holds for all $f \in H^{m+1}(\Gamma){\cap \widetilde{H}^{s}(\Gamma)}$.
\end{lemma}
The second estimate for $\partial\Gamma \neq \emptyset$ follows by extending $\varPi_hf \in V_h^p$ by zero outside $\Gamma$ which allows to estimate a $\widetilde{H}^{\pm \frac{1}{2}}$ norm on the left hand side by standard Sobolev norms (see \cite{screen}).\\
Combining $\varPi_h$ and $\varPi_{\Delta t}$ one obtains as in Proposition 3.54 of \cite{Glaefke}:
\begin{lemma}\label{interp}
Let $f \in H^{s}_\sigma(\mathbb{R}^+, H^m(\Gamma){\cap \widetilde{H}^{r}(\Gamma)})$, $0<m\leq q+1$, $0<s\leq p+1$, $r\leq s$, $|l|\leq \frac{1}{2}$ such that $lr\geq 0$. Then if  $ l,r\leq 0$
\begin{align*}\label{eq:approx}
\|f-\varPi_h \circ \varPi_{\Delta t} f\|_{r,l,\Gamma} &\leq C (h^\alpha + (\Delta t)^\beta)||f||_{s,m,\Gamma}\ ,\\
\|f-\varPi_h \circ \varPi_{\Delta t} f\|_{r,l,\Gamma,\ast} &\leq C \log_\ast(h) (h^\alpha + (\Delta t)^\beta)||f||_{s,m,\Gamma}\ ,
\end{align*}
where $\alpha = \min\{m-l, m-\frac{m(l+r)}{m+s}\}$, $\beta = \min\{m+s-(l+r), m+s-\frac{m+s}{m}l\}$. If $l,r>0$, $\beta = m+s-(l+r)$.
\end{lemma}
We are also going to require inverse estimates like (3.182) in \cite{Glaefke} for $s,m \leq 0$
$$\|p_{h,\Delta t}\|_{0,0,\Gamma}\leq C (\Delta t)^s \max{(h^m , \Delta t^m)} \|p_{h,\Delta t}\|_{s,m,\Gamma}\ $$
for $p_{h,\Delta t}$ in the approximation spaces $V_{h,\Delta t}^{p,q}$, namely
\begin{equation*}
 ||p_{h,\Delta t}||_{1,-\frac{1}{2},\Gamma,\ast} \lesssim \frac{1}{\Delta t} ||p_{h,\Delta t}||_{0,-\frac{1}{2},\Gamma,\ast}\quad \text{ (in the proof of the Theorem 6.1)}
\end{equation*}
\begin{equation*}
 ||p_{h,\Delta t}||_{1,0,\Gamma} \lesssim \frac{1}{\Delta t} ||p_{h,\Delta t}||_{0,0,\Gamma}\quad \text{ (in the proof of the Theorem 6.2)}
\end{equation*}
\begin{equation*}
 ||p_{h,\Delta t}||_{0,\frac{1}{2},\Gamma} \lesssim \frac{1}{\min\{\sqrt{\Delta t}, \sqrt{h}\}} ||p_{h,\Delta t}||_{0,0,\Gamma}\quad \text{ (in the proof of the Theorem 6.2)}.
\end{equation*}
The above inverse inequalities hold due to the standard estimates for regular finite element functions in the usual Sobolev spaces $H^s(\Gamma)$ \cite{Babuska_aziz} on one hand, and on the other hand the weight function $e^{-\sigma t}$ does not affect these inequalities (see \cite[Lemma 2]{bamberger86} ).

\section{Frequency--domain integral operators in the absorbing half--space}

We follow the approach by Bamberger and Ha-Duong \cite{bamberger86} and first analyze an associated Helmholtz problem in the frequency domain. The analysis will be translated into results for the wave equation in the following section.\\

Let $\sigma>0$. For a fixed frequency $\omega$ with $\mathrm{Im}\ \omega >\sigma$ we consider the exterior Helmholtz problem associated to the wave equation \eqref{eq:strong_half_space} for $u^e \in H^1(\Omega^e)$:
\begin{equation}\label{eq:exterior_helmholz_half}
\begin{cases}
(\Delta +\omega^2) u^e(x)=0 \quad\mbox{ in }\Omega^e \\
\frac{\partial u^e}{\partial n} +\alpha i \omega  u^e = \tilde{f}\quad\mbox{ on }\Gamma \\
\frac{\partial u^e}{\partial n} +\alpha_\infty i \omega  u^e  = 0\quad\mbox{ on } \Gamma_\infty
\end{cases}
\end{equation}
plus a Sommerfeld radiation condition at infinity.
The radiation condition holds automatically since for $\mathrm{Im}\ \omega >\sigma$ the solution decays like $e^{-\sigma|x|}$, and hence the solution belongs to $H^1(\Omega^e)$ and not only $H^1_{loc}(\Omega^e)$.\\
We also need an auxiliary interior problem for a function $u^i \in H^1(\Omega^i)$:
\begin{equation}\label{eq:interior_helmholz_half}
\begin{cases}
(\Delta +\omega^2) u^i(x)=0 \quad\mbox{ in }\Omega^i\\
\frac{\partial u^i}{\partial n} -\alpha i \omega  u^i = \tilde{g}\quad\mbox{ on }\Gamma\\
\frac{\partial u^i}{\partial n} +\alpha_\infty i \omega  u^i  = 0\quad\mbox{ on } \Gamma_\infty'= \partial\mathbb{R}^d_+ \setminus \Gamma_\infty \;.
\end{cases}
\end{equation}
The right-hand sides $\tilde{f}, \tilde{g}$ belong to ${H}^{-\frac{1}{2}}(\Gamma)$. \\

In the appendix we prove a uniqueness result:
\begin{theorem}
The problems \eqref{eq:exterior_helmholz_half}-\eqref{eq:interior_helmholz_half} admit at most one solution for $\mathrm{Re} \ \alpha \geq 0$ and $\mathrm{Re} \ \alpha_\infty \geq 0$.
\end{theorem}
The next step is to explicitly construct and represent the solution of the Helmholtz equation in $\Omega^e$ and $\Omega^i$
by means of layer potentials using the representation formula.\\

As derived by Ochmann \cite{Ochmann02}, for $d=3$ a fundamental solution to the half--space problem is given by:
$$G_\omega(x,y) = \frac{e^{i\omega |x-y|}}{4\pi |x-y|} + \frac{e^{i\omega |x-y'|}}{4\pi |x-y'|} + 2 \beta_\infty e^{-\beta_\infty (x_3+y_3)} \int_{\infty}^{-(x_3+y_3)} e^{-\beta_\infty \eta } \frac{e^{ik r(\eta)}}{4\pi r(\eta)} d\eta\ ,$$
where $r(\eta)=\sqrt{(x_1-y_1)^2+(x_2-y_2)^2+\eta^2}$ and $\beta_\infty=i\omega \alpha_\infty$. {For $y=(y_1,y_2,y_3)\in \mathbb{R}^3_+$, $y'$ is given by $y'=(y_1,y_2,-y_3)$.}
In any dimension, $G_\omega$ allows to define the potential operators for the absorbing half--space as
\begin{align*}
 S_\omega p (x)= \int_{\Gamma}G_\omega(x,y)\ p(y) \, ds_y\ ,\  D_\omega \varphi(x)=\int_{\Gamma}\frac{\partial G_\omega}{\partial n_y}(x,y)\ \varphi(y) \, ds_y\, .
\end{align*}

Using $G_\omega$, the solution $u$ of the Helmholtz problems admits
 an integral representation formula over $\Gamma$, not just $\Gamma\cup \Gamma_\infty\cup \Gamma_\infty'$.
\begin{theorem}
Any solution $u \in H^1(\Omega^i) \cup H^1(\Omega^e)$ of \eqref{eq:exterior_helmholz_half}-\eqref{eq:interior_helmholz_half} satisfying the acoustic boundary conditions on $\Gamma_\infty \cup \Gamma_\infty'$ admits a representation
\begin{equation*}
 u=S_\omega p- D_\omega \varphi  \quad \text{  in  } \quad \Omega^i \cup \Omega^e\, ,
\end{equation*}
where
\begin{equation*}
\varphi=u^i-u^e \quad \text{  and  } \quad p=\frac{\partial u^i}{\partial n} -\frac{\partial u^e}{\partial n} \quad \text{  on  } \Gamma\,.
\end{equation*}
\end{theorem}
The proof of the representation formula is standard if $\Gamma$ is replaced by $\Gamma\cup\Gamma_\infty\cup \Gamma_\infty'$ in the definition of $S_\omega$ and $D_\omega$. The contribution from integrals over $\Gamma_\infty$ and $\Gamma_\infty'$, however, vanishes since $G_\omega$ satisfies the acoustic boundary conditions.\\

Taking boundary values of $S_\omega$ and $D_\omega$, we obtain integral operators on $\Gamma$,
\begin{align*}
 & V_\omega p (x)= 2 \int_{\Gamma} G_\omega(x,y)\ p(y) \, ds_y \ , \quad & K'_\omega \varphi(x) = 2\int_{\Gamma} \frac{\partial}{\partial n_x} G_\omega(x,y)\ \varphi(y) \, ds_y \ ,\\
& K_\omega p(x)=2 \int_{\Gamma}\frac{\partial}{\partial n_y}G_\omega(x,y)\ p(y) \, ds_y \ , & W_\omega \varphi(x)= 2\int_{\Gamma}\frac{\partial^2}{\partial n_x \partial n_y} G_\omega(x,y)\ \varphi(y)  \, ds_y \ .
\end{align*}
{Here and in the following, the integrals are interpreted as distributional pairings, equivalently as principal values.} As in the full space, the operators relate the traces of $u$ with $\varphi$ and $p$:
\begin{align}\label{eq:trace_equation_freq}
 &2u^e =V_\omega p-(I +K_\omega)\varphi \ ,  &2u^i =V_\omega p+(I-K_\omega)\varphi \\
&2\frac{\partial u^e}{\partial n} = (-I + K'_\omega)p-W_\omega \varphi \ ,  &2\frac{\partial u^i}{\partial n} = (I + K'_\omega)p-W_\omega \varphi \,.\nonumber
\end{align}
Adding and subtracting the boundary conditions \eqref{eq:exterior_helmholz_half}-\eqref{eq:interior_helmholz_half} on $\Gamma$, we have
\begin{equation*}
\begin{cases}
\frac{\partial u^e}{\partial n} + \frac{\partial u^i}{\partial n}-\alpha i \omega \varphi=\tilde{f}+\tilde{g}=F\\
p-\alpha i \omega(u^e+u^i)=\tilde{g}-\tilde{f}=G\,.
\end{cases}
\end{equation*}
Then using the equation \eqref{eq:trace_equation_freq} of the trace $u$ we find the following system of  integral equations:
\begin{equation}\label{eq:sys_op_freq}
\begin{cases}
 K'_\omega p -W_\omega \varphi -i\omega \alpha  \varphi =F\\
p- i \omega \alpha (V_\omega p -K_\omega \varphi)=G\,.
\end{cases}
\end{equation}
If $\alpha\neq 0$ multiplying the first equation by $\overline{-i\omega \psi}$ and the second by $\frac{1}{\alpha} \bar{q}$, we obtain the weak formulation after an integration by parts: Find $\tilde{\Phi} = (\varphi, p)$ such that
\begin{equation}\label{eq:var_freq}
 a_\omega(\tilde{\Phi},\tilde{\Psi})=l_\omega (\tilde{\Psi}) \qquad \text{for all $\tilde{\Psi} = (\psi, q)$}.
\end{equation}
Here,
\begin{align*}
a_\omega(\tilde{\Phi},\tilde{\Psi})&=|\omega|^2 \int_\Gamma \alpha \varphi \bar{\psi} ds_x + \int_\Gamma \frac{1}{\alpha} p \bar{q} ds_x +  i \bar{\omega}\int_\Gamma K'_\omega p \bar{\psi} ds_x\\
&\qquad - i\bar{\omega}\int_\Gamma W_\omega \varphi \bar{\psi} ds_x - i\omega \int_\Gamma V_\omega p \bar{q} ds_x + i\omega \int_\Gamma K_\omega \varphi \bar{q} ds_x
\end{align*}
and $l_\omega(\tilde{\Psi})= i\bar{\omega} \int_\Gamma F \bar{\psi} ds_x +\int_\Gamma  \frac{1}{\alpha} G\bar{q} ds_x$.
For $\alpha= 0$, \eqref{eq:sys_op_freq} reduces to $W_\omega \varphi=K'_\omega G-F $. {For simplicity, we assume $\alpha^{-1}$ to exist. Other cases have to be treated differently.}
\begin{theorem}\label{3.3}
 \emph{(Coercivity)}\\
Assume that $\mathrm{Re}\ \alpha>0,\mathrm{Re}\ \alpha_\infty\geq0$. Then the following inequality holds for all $\tilde{U}=(\varphi,p) \in \widetilde{H}^{\frac{1}{2}}(\Gamma)\times L^2(\Gamma)$:
\begin{equation*}
 \mathrm{Re}\ a_\omega(\tilde{U},\tilde{U}) \gtrsim_\sigma \|(\mathrm{Re}\ \alpha^{-1})^{1/2} p\|^2_{0,\omega, \Gamma}+ ||\varphi||^2_{\frac{1}{2},\omega,\Gamma,\ast}+  ||\omega (\mathrm{Re}\ \alpha)^{1/2}\varphi||^2_{0,\omega, \Gamma}\,.
\end{equation*}
\end{theorem}
\noindent \emph{Proof: }
Taking the real part of the bilinear form $a_\omega$ and using \eqref{eq:trace_equation_freq}, we calculate
\begin{equation*}
 \mathrm{Re}(a_\omega(\tilde{U},\tilde{U}))=\mathrm{Re} \int_\Gamma (K'_\omega p -W_\omega \varphi-i\omega \alpha  \varphi)(\overline{-i\omega \varphi}) + \bar{p} \frac{p-i \omega \alpha(V_\omega p -K_\omega \varphi)}{\alpha}\,ds_x
\end{equation*}
\begin{align*}
=\mathrm{Re} \int_\Gamma [\frac{\partial u^i}{\partial n} &+ \frac{\partial u^e}{\partial n}-i \omega \alpha(u^i-u^e)]i\bar{\omega} (\bar{u^i}-\bar{u^e}) \,ds_x\\
          &+\mathrm{Re} \int_\Gamma \frac{1}{\alpha} (\frac{\partial \bar{u^i}}{\partial n} - \frac{\partial \bar{u^e}}{\partial n})(\frac{\partial u^i}{\partial n} - \frac{\partial u^e}{\partial n}-i\omega \alpha (u^i+u^e))\,ds_x\\
\quad =\mathrm{Re} \int_\Gamma \,i\bar{\omega}(2\frac{\partial u^i}{\partial n} \bar{u^i} - &2\frac{\partial u^e}{\partial n}\bar{u^e})\,ds_x\\
          &+ \int_{\Gamma} \frac{1}{\alpha} \underbrace{|\frac{\partial \bar{u^i}}{\partial n} - \frac{\partial \bar{u^e}}{\partial n}|^2}_{=|p|^2} \,ds_x+ |\omega|^2 \int_{\Gamma} \alpha \underbrace{|u^i-u^e|^2}_{=|\varphi|^2}\,ds_x\,.
\end{align*}
{Adding $0$},
\begin{align*}
 -\int_\Gamma \frac{\partial u^e}{\partial n}\bar{u^e} \,ds_x&=-\int_\Gamma \frac{\partial u^e}{\partial n}\bar{u^e} \,ds_x -\int_{\Gamma_\infty} \frac{\partial u^e}{\partial n}\bar{u^e} \,ds_x +\int_{\Gamma_\infty}{\frac{\partial u^e}{\partial n}}\bar{u^e} \,ds_x\,,
\end{align*}
integration by parts on $\Omega^e$ leads to:
\begin{align*}
 -\int_\Gamma \frac{\partial u^e}{\partial n}\bar{u^e} \,ds_x &=\int_{\Omega^e} \Delta u^e \,\bar{u^e}+ \triangledown u^e \,\overline{\triangledown u^e} \,dx +\int_{\Gamma_\infty} \frac{\partial u^e}{\partial x_3} \bar{u^e}\,ds_x\\
 &=\int_{\Omega^e} |\triangledown u^e|^2-\omega^2 |u^e|^2 \,dx + \int_{\Gamma_\infty} \frac{\partial u^e}{\partial x_3} \bar{u^e}\,ds_x\\
&=\int_{\Omega^e} |\triangledown u^e|^2-\omega^2 |u^e|^2  \,dx- \int_{\Gamma_\infty} i\alpha_\infty \omega |u^e|^2 \,ds_x\,.
\end{align*}
Therefore,
\begin{align*}
 -\mathrm{Re}\, 2 i \bar{\omega} \int_\Gamma \frac{\partial u^e}{\partial n}\bar{u^e} \,ds_x &= \mathrm{Re} (2 \int_{\Omega^e} i\bar{\omega}|\triangledown u^e|^2- i \bar{\omega} \omega^2 |u^e|^2 \,dx -2\int_{\Gamma_\infty} (i\bar{\omega}) i\alpha_\infty \omega |u^e|^2 \,ds_x)\\
&=\mathrm{Re} (2 \int_{\Omega^e} i\bar{\omega}|\triangledown u^e|^2-i \omega |\omega|^2 |u^e|^2  \,dx +2 \int_{\Gamma_\infty} \alpha_\infty |\omega|^2 |u^e|^2 \,ds_x)\\
&{\geq}2 \sigma \int_{\Omega^e} |\triangledown u^e|^2+|\omega|^2 |u^e|^2 \,dx +2 (\mathrm{Re}\ \alpha_\infty) \int_{\Gamma_\infty} |\omega|^2 |u^e|^2 \,ds_x\,.
\end{align*}
Similarly,
\begin{align*}
 \mathrm{Re}\, 2 i \bar{\omega} \int_\Gamma \frac{\partial u^i}{\partial n}\bar{u^i} \,ds_x &= \mathrm{Re} (2 \int_{\Omega^i} i\bar{\omega}|\triangledown u^i|^2-i \bar{\omega}\omega^2 |u^i|^2 \,dx +2 \int_{\Gamma_\infty'} \alpha_\infty |\omega|^2 |u^i|^2 \,ds_x)\\
&{\geq}2 \sigma \int_{\Omega^i} |\triangledown u^i|^2+|\omega|^2 |u^i|^2 \,dx +2 (\mathrm{Re}\ \alpha_\infty) \int_{\Gamma_\infty'} |\omega|^2 |u^i|^2 \,ds_x\,.
\end{align*}
We conclude
\begin{align*}
\mathrm{Re}\ a_\omega(\tilde{U},\tilde{U})&=\mathrm{Re} \,2 i\bar{\omega}\int_\Gamma (\frac{\partial u^i}{\partial n} \bar{u^i} - \frac{\partial u^e}{\partial n}\bar{u^e}) \,ds_x+   \int_{\Gamma} \frac{1}{\alpha}|p|^2 +  |\omega|^2 \int_{\Gamma}\alpha |\varphi|^2 \,ds_x\\
&{\geq}2 \sigma \int_{\Omega^i\cup \Omega^e} |\triangledown u|^2+|\omega|^2 |u|^2\,dx + \int_{\Gamma} \mathrm{Re}(\frac{1}{\alpha})|p|^2 \,ds_x +  |\omega|^2 \int_{\Gamma}\mathrm{Re}(\alpha) |\varphi|^2 \,ds_x\\
&\qquad\qquad \qquad \qquad +2\int_{\Gamma_\infty} \mathrm{Re}(\alpha_\infty) |\omega|^2 |u^e|^2 \,ds_x+2 \int_{\Gamma_\infty'} \mathrm{Re}(\alpha_\infty) |\omega|^2 |u^i|^2 \,ds_x\\
&\geq 2 \sigma \int_{\Omega^i\cup \Omega^e} |\triangledown u|^2+|\omega|^2 |u|^2 \,dx + \int_{\Gamma}\mathrm{Re}(\frac{1}{\alpha}) |p|^2 \,ds_x +  |\omega|^2 \int_{\Gamma} \mathrm{Re}(\alpha) |\varphi|^2 \,ds_x\,.
\end{align*}
Using the trace theorem in $\Omega^i$ and $\Omega^e$, $||\varphi||_{\frac{1}{2},\omega,\Gamma,\ast}\lesssim_\sigma ||u||_{1,\omega,\Omega}$, we obtain the assertion:
\begin{equation*}
\mathrm{Re}\ a_\omega(\tilde{U},\tilde{U}) \gtrsim_\sigma \|(\mathrm{Re}\ \alpha^{-1})^{1/2} p\|^2_{0,\omega, \Gamma}+ ||\varphi||^2_{\frac{1}{2},\omega,\Gamma,\ast}+  ||\omega (\mathrm{Re}\ \alpha)^{1/2}\varphi||^2_{0,\omega, \Gamma}\,.
\end{equation*}

\begin{remark}
 Assume $\mathrm{Re} \ \alpha_\infty \geq 0$. Then a similar coercivity estimate holds for the single layer potential $V_\omega$:
\begin{equation}\label{eq:v_omega_coerc}
  \mathrm{Re} \langle i\omega V_\omega {\varphi}, \varphi\rangle \geq C_\sigma ||\varphi||^2_{-\frac{1}{2},\omega,\Gamma,\ast}\,.
\end{equation}
\end{remark}
Boundedness of the integral operators is also shown by going into $\Omega^e\cup \Omega^i$. We postpone the proof to the appendix.
\begin{theorem}\label{3.4}
\emph{(Continuity)}\\
Assume that $\mathrm{Re}\ \alpha_\infty\geq0$. The integral operators satisfy the following mapping properties for $p\in \widetilde{H}^{-\frac{1}{2}}(\Gamma)$ and $\varphi \in \widetilde{H}^{\frac{1}{2}}(\Gamma)$:
 \begin{equation}\label{eq:cont_V}
 ||V_\omega p||_{\frac{1}{2},\omega,\Gamma} \lesssim_\sigma |\omega|||p||_{-\frac{1}{2},\omega,\Gamma,\ast}\,,
\end{equation}
\begin{equation}\label{eq:cont_W}
 ||W_\omega \varphi||_{-\frac{1}{2},\omega,\Gamma} \lesssim_\sigma |\omega|||\varphi||_{\frac{1}{2},\omega,\Gamma,\ast}\,,
\end{equation}
\begin{equation}\label{eq:cont_K}
 ||(I-K_\omega)\varphi||_{\frac{1}{2},\omega,\Gamma} \lesssim_\sigma |\omega|||\varphi||_{\frac{1}{2},\omega,\Gamma,\ast}\,,
\end{equation}
\begin{equation}\label{eq:cont_K'}
 ||(I-K'_\omega)p||_{-\frac{1}{2},\omega,\Gamma} \lesssim_\sigma |\omega|||p||_{-\frac{1}{2},\omega,\Gamma,\ast}\,.
\end{equation}
\end{theorem}
The theorem translates into the boundedness of the considered bilinear form.
\begin{theorem}\label{3.5}
 Assume that $\mathrm{Re}\ \alpha_\infty\geq 0$ and $\alpha,\, \frac{1}{\alpha} \in L^\infty(\Gamma)$.
 The bilinear form $a_\omega$ is continuous on $\left(\widetilde{H}^{\frac{1}{2}}(\Gamma)\times L^2(\Gamma)\right)\times \left(\widetilde{H}^{\frac{1}{2}}(\Gamma)\times L^2(\Gamma)\right)$.
\end{theorem}

Again, we refer to the appendix for a proof.\\

{With these results we can now state the precise weak formulation of the boundary integral equation \eqref{eq:var_freq}:
Find $\tilde{\Phi}=(\varphi,p)\in \widetilde{H}^{\frac{1}{2}}(\Gamma) \times L^2(\Gamma)$ such that for all $\tilde{\Psi}=(\psi,q)\in \widetilde{H}^{\frac{1}{2}}(\Gamma) \times L^2(\Gamma)$:
\begin{equation*}
 a_\omega(\tilde{\Phi},\tilde{\Psi})=l_\omega (\tilde{\Psi}).
\end{equation*}
Using the coercivity estimate in Theorem \ref{3.3}, we conclude the following estimate on the solution:
\begin{equation}\label{omegaacousticbound}\|p\|_{0,\omega, \Gamma}+ ||\varphi||_{\frac{1}{2},\omega,\Gamma,\ast}+  ||\omega \varphi||_{0,\omega, \Gamma} \lesssim_\sigma \min\{||\omega F||_{-\frac{1}{2},\omega,\Gamma}, ||F||_{0,\omega,\Gamma}\} + ||G||_{0,\omega,\Gamma}\ .
\end{equation}
Similar results are obtained for the weak formulation of the Dirichlet problem, which reads: Find $\phi \in \widetilde{H}^{-\frac{1}{2}}(\Gamma)$ such that for all $\psi \in \widetilde{H}^{-\frac{1}{2}}(\Gamma)$:
\begin{equation*}
\langle V_\omega\phi,\psi\rangle=\langle f,\psi\rangle .
\end{equation*}
When $\mathrm{Re} \ \alpha_\infty \geq 0$, from the coercivity \eqref{eq:v_omega_coerc} one obtains the estimate
 \begin{equation}\label{omegadirbound}
 ||\phi||_{-\frac{1}{2},\omega,\Gamma, \ast} \lesssim_\sigma |\omega|||f||_{\frac{1}{2},\omega,\Gamma}
\end{equation}
on the solution.
}

\section{Time--domain boundary integral equations for an absorbing half--space}

We consider the wave equation in $\mathbb{R}^d_+$ with acoustic boundary condition
\begin{equation*}
 \frac{\partial u}{\partial n}-\alpha_\infty \frac{\partial u}{\partial t}=0\, \qquad \text{ on $\partial\mathbb{R}^d_+$}\ .
\end{equation*}
In $\mathbb{R}^3_+$ Ochmann determines the Green's function to be \cite{Ochmann01}
\begin{align}\label{eq:green_absorbing}
 G(t-s,x,y)=\frac{\delta(t-s-r(y_3))}{4\pi r(y_3)}+\frac{\delta(t-s-r(-y_3))}{4\pi r(-y_3)} + \Sigma
\end{align}
with
\begin{align*}
 \Sigma=\frac{-\alpha_\infty}{2\pi} \frac{\partial}{\partial t} \frac{H(t-s-r(-y_3))}{\sqrt{(t-s +\alpha_\infty(x_3+y_3))^2+(\alpha_\infty^2-1)R^2}}\ .
\end{align*}
Here $H$ denotes the Heaviside function, $R^2=(x_1-y_1)^2+(x_2-y_2)^2$ and $r(\pm y_3)^2 = R^2+(x_3\mp y_3)^2$. The second and third terms on the right-hand side of ${G}$ represent the field reflected by the plane $\Gamma_\infty$. After a Fourier transform in $t$, one recovers from $G$ the frequency--domain Green's function $G_\omega$ from Section 4.\\

As for the Helmholtz problem the solution $u$ of the direct scattering problem \eqref{eq:strong_half_space} and its associated interior problem admits
 an integral representation formula over $\Gamma$, not just $\Gamma\cup \Gamma_\infty\cup\Gamma_\infty'$. A similar representation formula in time--domain has been obtained by Becache \cite{becache94} for exterior domains in $\mathbb{R}^3$.
\begin{theorem}
Let $u \in L^2(\mathbb{R}^+,H^1(\Omega^i \cup \Omega^e))\cap H^1_0(\mathbb{R}^+, L^2(\Omega^i \cup \Omega^e))$ be the solution of \eqref{eq:strong_half_space} for a Lipschitz boundary $\Gamma$. Then it holds in the sense of distributions ($x \in \Omega^e\cup \Omega^i$, $t \in \mathbb{R}^+$):
\begin{align*}
u(t,x) &=  \int_{\mathbb{R}^+ \times \Gamma} \frac{\partial G}{\partial n_{y}}(t- \tau,x,y) u(\tau,y) d\tau ds_y  \\
        &\quad   - \int_{\mathbb{R}^+ \times \Gamma} G(t- \tau,x,y) \frac{\partial u}{\partial n_y} (\tau,y) d\tau ds_y \ ,
\end{align*}
where $G$ is a fundamental solution in the half-space which satisfies the acoustic boundary conditions.
\end{theorem}

We introduce the single layer potential in time domain for a half-space with an absorbing boundary condition as
$$ S p(t,x)=\int_{\mathbb{R}^+ \times \Gamma} G(t- \tau,x,y) p(\tau,y) d\tau ds_y\ .$$
Specifically in 3 dimensions, this is
\begin{align*}
S p(t,x) &=\frac{1}{4 \pi} \int_\Gamma \frac{p(t-|x-y|,y)}{|x-y|} ds_y + \frac{1}{4 \pi} \int_\Gamma  \frac{p(t-|x-y'|,y)}{|x-y'|} ds_y \\
         &\qquad -\frac{\alpha_\infty}{2 \pi} \int\limits_0^\infty\int_\Gamma \frac{\partial}{\partial s}  \Big[\frac{H(t-s-|x-y'|)}{\sqrt{(t-s +\alpha_\infty (x_3+y_3))^2+(\alpha_\infty^2-1)R^2}}\Big] p(s,y) ds_y ds\,.
\end{align*}
The corresponding double layer potential $D$ is:
\begin{align*}
 D \varphi(t,x)&=\int_{\mathbb{R}^+ \times \Gamma} \frac{\partial G}{\partial n_{y}}(t- \tau,x,y)\varphi(\tau,y) d\tau ds_y\,.
\end{align*}
The function $u=Sp-D\varphi$ satisfies the wave equation on $\mathbb{R}^d_+ \setminus \Gamma$, and according to the representation formula
\begin{equation*}
 \varphi= u^i-u^e\,,\quad p=\frac{\partial u^i}{\partial n}-\frac{\partial u^e}{\partial n}  \quad\mbox{ on  } \mathbb{R}^+ \times \Gamma\,.
\end{equation*}
As for the Helmholtz equation we have the following trace identities:
\begin{align}\label{eq:rep_trace}
& 2u^e =V p -(I+K)\varphi \qquad & 2u^i =V p+(I-K)\varphi \ , \\
&2\frac{\partial u^e}{\partial n} =(-I+K')p- W \varphi &2\frac{\partial u^i}{\partial n} = (I+K')p- W \varphi \,.\notag
\end{align}
The relevant boundary integral operators on $\Gamma$ are:
\begin{align*}
& V p (t,x)= 2 \int_{\mathbb{R}^+ \times \Gamma} G(t- \tau,x,y) p(\tau,y) d\tau ds_y\, ,\\ &K' \varphi(t,x)=2 \int_{\mathbb{R}^+\times \Gamma} \frac{\partial G}{\partial n_x}(t- \tau,x,y)\varphi(\tau,y) d\tau ds_y\, ,\\
&K\varphi(t,x)=2 \int_{\mathbb{R}^+\times \Gamma} \frac{\partial G}{\partial n_y}(t- \tau,x,y)\varphi(\tau,y) d\tau ds_y,\ \\ &W \varphi(t,x)=2 \int_{\mathbb{R}^+\times \Gamma} \frac{\partial^2 G}{\partial n_x \partial n_y}(t- \tau,x,y)\varphi(\tau,y) d\tau ds_y \ .
\end{align*}
{Because $G$ is the inverse Fourier--Laplace transform $\mathcal{F}^{-1}_{\omega \to t} G_\omega$, see \cite{Ochmann01}, these boundary integral operators are conjugates of their frequency--domain analogues: $V = \mathcal{F}^{-1}_{\omega \to t}\circ V_\omega \circ\mathcal{F}_{t \to \omega}$, and analogously for $K', K, W$.\\} Substituting formula \eqref{eq:rep_trace} into the boundary condition on $\Gamma$, we obtain the following system for the unknown functions $\varphi$ and $p$
\begin{equation}\label{eq:sys_op}
\begin{cases}
(-I+K')p-W \varphi-\alpha \partial_t(Vp-(I+K)\varphi)=2f \\
(I+K')p-W \varphi+\alpha \partial_t(Vp+(I-K)\varphi)=2g \,.
\end{cases}
\end{equation}
Adding respectively subtracting the two equations of \eqref{eq:sys_op}, again leads to
\begin{equation}\label{eq:sys_op_time}
\begin{cases}
K' p -W \varphi + \alpha \frac {\partial \varphi}{\partial t}=F\\
p+ \alpha (V \partial_t p-K {\partial_t \varphi})=G\,.
\end{cases}
\end{equation}
Pairing these equations with test functions $\partial_t\psi$ respectively $\frac{q}{\alpha}$, we obtain the following space-time variational formulation:
\begin{align*}
 \int_0^\infty \int_\Gamma \left[(K' p -W \varphi) + \alpha \partial_t{\varphi}\right] \partial_t{\psi}\, ds_x\, d_\sigma t &= \int_0^\infty \int_\Gamma F \partial_t{\psi} \,ds_x\, d_\sigma t \\
 \int_0^\infty \int_\Gamma \left[\frac{p}{\alpha}+ (V \partial_t p-K  \partial_t \varphi)\right] q \,ds_x\, d_\sigma t &= \int_0^\infty \int_\Gamma \frac{G q}{\alpha}\, ds_x\, d_\sigma t\,.
\end{align*}
Here $d_\sigma t = e^{-2\sigma t} dt$, $\sigma>0$. The system can be written as
\begin{equation}\label{eq:var_time_space}
 a(\Phi,\Psi)=l(\Psi) \,,
\end{equation}
where  $\Phi=(\varphi,p)$, $\Psi=(\psi,q)$  and
\begin{equation}\label{eq:bilinear_time}
 a(\Phi,\Psi)=\int_0^\infty \int_\Gamma \left( \alpha (\partial_t{\varphi}) (\partial_t{\psi})  +  \frac{1}{\alpha} p q  +  K' p (\partial_t{\psi})
- W \varphi (\partial_t\psi) + V (\partial_t p) q -K (\partial_t\varphi) q \right) ds_x\, d_\sigma t \,,
\end{equation}
\begin{equation}\label{eq:rechte_time}
 l(\Psi)=\int_0^\infty \int_\Gamma  F \partial_t\psi \,ds_x\, d_\sigma t+\int_0^\infty \int_\Gamma \frac{G q}{\alpha}\, ds_x\, d_\sigma t\,.
\end{equation}
\begin{remark}
 The system of equations \eqref{eq:sys_op_time} and the variational formulation \eqref{eq:var_time_space} are the inverse Fourier-Laplace  transforms of \eqref{eq:sys_op_freq} and \eqref{eq:var_freq}.
\end{remark}

Later we will also require the time--domain mapping properties of the boundary integral operators in the energy Sobolev spaces.

\begin{theorem}\label{mapthm}
The following operators are continuous for $r \in \mathbb{R}$:
\begin{align*}
V&: H_\sigma^{r+1}(\mathbb{R}^+, \widetilde{H}^{-\frac{1}{2}}(\Gamma)) \to H_\sigma^{r}(\mathbb{R}^+, H^{\frac{1}{2}}(\Gamma)) \ ,\\ K' &: H_\sigma^{r+1}(\mathbb{R}^+, \widetilde{H}^{-\frac{1}{2}}(\Gamma)) \to H_\sigma^{r}(\mathbb{R}^+, H^{-\frac{1}{2}}(\Gamma))\ ,\\
K&: H_\sigma^{r+1}(\mathbb{R}^+, \widetilde{H}^{\frac{1}{2}}(\Gamma)) \to H_\sigma^{r}(\mathbb{R}^+, H^{\frac{1}{2}}(\Gamma)) \ ,\\ W &: H_\sigma^{r+1}(\mathbb{R}^+, \widetilde{H}^{\frac{1}{2}}(\Gamma)) \to H_\sigma^{r}(\mathbb{R}^+, H^{-\frac{1}{2}}(\Gamma))\ .
\end{align*}
\end{theorem}
\noindent \emph{Proof: }
{Like the corresponding assertions in the full space \cite{Ha-Duong03a},} the theorem follows from Theorem {\ref{3.4}} {and Definition \ref{sobdef} by conjugation with the Fourier transform: $V = \mathcal{F}^{-1}_{\omega \to t}\circ V_\omega \circ\mathcal{F}_{t \to \omega}$, and analogously for $K', K, W$.\\}.\\

Together with Theorem \ref{3.3}, the mapping properties imply continuity and coercivity of the bilinear form $a(U,V)$.
\begin{theorem}\label{conti_coer_time}
 Assume that $\mathrm{Re}\ \alpha_\infty\geq 0$ and $\alpha,\,\frac{1}{\alpha} \in L^\infty(\Gamma)$. Then the bilinear form of the variational formulation \eqref{eq:var_time_space} is continuous on $\left(H_\sigma^1(\mathbb{R}^+,\widetilde{H}^{\frac{1}{2}}(\Gamma)) \times H_\sigma^1(\mathbb{R}^+,L^2(\Gamma))\right) \times \left(H_\sigma^1(\mathbb{R}^+,\widetilde{H}^{\frac{1}{2}}(\Gamma)) \times H_\sigma^1(\mathbb{R}^+,L^2(\Gamma))\right) $, i.e.
\begin{equation}\label{eq:var_cont}
 |a(U,V)|\lesssim_\sigma(||p||_{1,0,\Gamma}+||\varphi||_{1,\frac{1}{2},\Gamma,\ast})(||q||_{1,0,\Gamma}+||\psi||_{1,\frac{1}{2},\Gamma,\ast}) \,.
\end{equation}
If $\mathrm{Re}\ \alpha>0$, it verifies a coercivity estimate: There exists $C_\sigma > 0$ such that:
\begin{equation}\label{eq:var_coerc}
 a(U,U) \geq C_\sigma (||p||^2_{0,0,\Gamma}+||\varphi||^2_{0,\frac{1}{2},\Gamma,\ast}+||\partial_t\varphi||^2_{0,0,\Gamma})\,.
\end{equation}
\end{theorem}
\noindent \emph{Proof: }
Equations \eqref{eq:var_cont} and \eqref{eq:var_coerc} follow from Theorem \ref{3.3} and Theorem \ref{mapthm}.\\
Concerning \eqref{eq:var_coerc} we note that
\begin{align*}
a(U,U)&=|a(U,U)|=|\int_{-\infty +i\sigma}^{\infty +i\sigma}   a_\omega(\tilde{U},\tilde{U}) d\omega|\\
      &\geq |\int_{-\infty +i\sigma}^{\infty +i\sigma}  \mathrm{Re}\ a_\omega(\tilde{U},\tilde{U}) d\omega |\\
      & \gtrsim ||p||^2_{0,0,\Gamma}+||\varphi||^2_{0,\frac{1}{2},\Gamma,\ast}+||\partial_t\varphi||^2_{0,0,\Gamma}\,.
\end{align*}
Similarly \eqref{eq:var_cont} is a consequence of \eqref{eq:a_freq_cont} and Cauchy-Schwarz.

\begin{remark}\label{Vcoerce}
 Similarly for the Dirichlet problem (see \cite[(54)]{Ha-Duong03a} and Corollary 3.50 in \cite{Glaefke} for the full space)
\begin{equation}
b(\varphi,\varphi)=\int_0^\infty\int_\Gamma V (\partial_t\varphi(t,x)) \varphi(t,x) d\sigma_x \, d_\sigma t \geq C_\sigma \|\varphi\|_{0,-\frac{1}{2},\Gamma,\ast}^2
\end{equation}
\end{remark}

\section{A priori estimates in the absorbing half--space}

In the following, we will restrict ourselves to a polyhedral surface $\Gamma$ ($d=3$) resp.~a polygonal curve ($d=2$), which satisfies the assumptions from Section 3. The error incurred by approximating a general smooth surface or curve by a polyhedral one has been studied by Nedelec \cite{nedelec}, and by Bamberger and Ha Duong in the context of the wave equation \cite{bamberger86}.

\subsection{Dirichlet problem}

We now use the approximation results of Section 3 to discuss the convergence of Galerkin approximations to the Dirichlet problem.
{As in frequency domain, the mapping properties of the integral operators allow us to state the precise weak formulation of the time--dependent boundary integral equation $V \phi = f$:}
Find $\phi \in H^1_{\sigma}(\mathbb{R}^+,\widetilde{H}^{-\frac{1}{2}}(\Gamma))$ such that
\begin{equation}\label{eq:symm}
 b(\phi,\psi)=\langle \partial_t f,\psi\rangle \qquad \forall \psi \in H^1_{\sigma}(\mathbb{R}^+,\widetilde{H}^{-\frac{1}{2}}(\Gamma))\,,
\end{equation}
where
\begin{align*}
 b(\phi,\psi)&=\int_0^\infty\int_\Gamma (V \partial_t\phi(t,x)) \psi(t,x) ds_x \, d_\sigma t\ ,\\
 \langle \partial_t f,\psi \rangle&= \int_0^\infty \int_\Gamma (\partial_t f(t,x)) \psi(t,x) ds_x \, d_\sigma t\,.
\end{align*}
{Similar to the estimate \eqref{omegadirbound} we obtain
\begin{equation}
||\phi||_{r, -\frac{1}{2}, \Gamma, \ast} \lesssim_\sigma ||f||_{r+1, \frac{1}{2}, \Gamma}
\end{equation}
for any $r\in \mathbb{R}$, provided that $\mathrm{Re} \ \alpha_\infty\geq 0$.
In particular, a solution $\phi \in H^1_{\sigma}(\mathbb{R}^+,\widetilde{H}^{-\frac{1}{2}}(\Gamma))$ only exists provided $f \in H^2_{\sigma}(\mathbb{R}^+,H^{\frac{1}{2}}(\Gamma))$. Coercivity assures that the solution is unique in this case.
}

The Galerkin formulation of \eqref{eq:symm} reads: Find $\phi_{h,\Delta t} \in V_{h,\Delta t}^{p,q}$ such that
\begin{equation}\label{eq:Galerkin_1}
 b(\phi_{h,\Delta t},\psi_{h,\Delta t})= \langle (\partial_t f)_{h,\Delta t},\psi_{h,\Delta t} \rangle \qquad \forall \psi_{h,\Delta t} \in V_{h,\Delta t}^{p,q}\ .
\end{equation}

{
For the solutions of the continuous and discrete problems we obtain the following a priori error estimate:
}
\begin{theorem}
Assume that $\mathrm{Re} \ \alpha_\infty\geq 0$. For the solutions $\phi\in H^1_{\sigma}(\mathbb{R}^+,\widetilde{H}^{-\frac{1}{2}}(\Gamma))$ of \eqref{eq:symm} and $\phi_{h,\Delta t} \in V_{h,\Delta t}^{p,q}$ of \eqref{eq:Galerkin_1} there holds:
\begin{align*}
 \|\phi-\phi_{h,\Delta t}\|_{0,-\frac{1}{2},\Gamma,\ast} &\lesssim||(\partial_t f)_{h,\Delta t}-{\partial_t f}||_{0,\frac{1}{2},\Gamma}\\
           &  \qquad \qquad + \inf_{\psi_{h,\Delta t} } \left\{  (1 + \frac{1}{\Delta t})\|{\phi}-{\psi}_{h,\Delta t}||_{0,-\frac{1}{2},\Gamma,\ast}+\frac{1}{\Delta t}\|\partial_t\phi-\partial_t\psi_{h,\Delta t}||_{0,-\frac{1}{2},\Gamma,\ast}\right\}\ .
\end{align*}
If in addition $\phi\in H^{s}_\sigma(\mathbb{R}^+,  {H}^m(\Gamma))$, then
\begin{align*}
 \|\phi-\phi_{h,\Delta t}\|_{0,-\frac{1}{2},\Gamma,\ast} &\lesssim ||(\partial_tf)_{h,\Delta t}-{\partial_t f}||_{0,\frac{1}{2},\Gamma}\\
&  \qquad \qquad+\left((h^{\alpha_1}+\Delta t^{\beta_1})(1+\frac{1}{\Delta t})+(h^{\alpha_2}+\Delta t^{\beta_2})\frac{1}{\Delta t}\right)||\phi||_{s,m,\Gamma}\ ,
\end{align*}
where
\begin{align*}
   \alpha_1 = \min\{m+\frac{1}{2}, m-\frac{m}{2(m+s)}\}&, \beta_1 = \min\{m+s+\frac{1}{2}, m+s+\frac{m+s}{2m}\} \,,\\
   \alpha_2 = \min\{m+\frac{1}{2}, m-\frac{m}{2(m+s-1)}\}&, \beta_2 = \min\{m+s-\frac{1}{2}, m+s-1+\frac{m+s-1}{2m}\}\,,
\end{align*}
and $\quad m\geq -\frac{1}{2},\,s\geq 0$.
\end{theorem}
\noindent \emph{Proof: }
We apply the coercivity from Remark \ref{Vcoerce} to $\phi_{h,\Delta t}-\psi_{h,\Delta t} \in H^1_\sigma(\mathbb{R}^+,\widetilde{H}^{-\frac{1}{2}}(\Gamma))$, $\psi_{h,\Delta t} \in V_{h,\Delta t}$:
\begin{align*}
 \|\phi_{h,\Delta t}-\psi_{h,\Delta t}\|_{0,-\frac{1}{2},\Gamma,\ast}^2\lesssim b({\phi}_{h,\Delta t}-{\phi},\phi_{h,\Delta t}-\psi_{h,\Delta t})+b({\phi}-{\psi}_{h,\Delta t},\phi_{h,\Delta t}-\psi_{h,\Delta t})\,.
\end{align*}
Continuity of the duality pairing is used to estimate the first term:
\begin{align*}
 b({\phi}_{h,\Delta t}-{\phi},\phi_{h,\Delta t}-\psi_{h,\Delta t})&=\int_0^{\infty} \int_\Gamma ((\partial_t{f})_{h,\Delta t}-\partial_t f)(\phi_{h,\Delta t}-\psi_{h,\Delta t})\, ds_x \,d_\sigma t\\
                                                                &\leq \|(\partial_t f)_{h,\Delta t}-\partial_t f\|_{0,\frac{1}{2},\Gamma} \| \phi_{h,\Delta t}-\psi_{h,\Delta t}\|_{0,-\frac{1}{2},\Gamma,\ast}\,.
\end{align*}
The mapping properties of $V$ from Theorem \ref{mapthm} bound the second term as follows:
\begin{align*}
 b({\phi}-{\psi}_{h,\Delta t},\phi_{h,\Delta t}-\psi_{h,\Delta t}) &\leq  \|V\partial_t(\phi-\psi_{h,\Delta t})\|_{-1,\frac{1}{2},\Gamma}\| \phi_{h,\Delta t}-\psi_{h,\Delta t}\|_{1,-\frac{1}{2},\Gamma,\ast}\\
& \lesssim  \|\phi-\psi_{h,\Delta t}\|_{1,-\frac{1}{2},\Gamma,\ast}\| \phi_{h,\Delta t}-\psi_{h,\Delta t}\|_{1,-\frac{1}{2},\Gamma,\ast}\ .
\end{align*}
The inverse inequality in the time variable leads to
\begin{equation*}
 b({\phi}-{\psi}_{h,\Delta t},\phi_{h,\Delta t}-\psi_{h,\Delta t}) \lesssim \frac{1}{\Delta t} \|\phi-\psi_{h,\Delta t}\|_{1,-\frac{1}{2},\Gamma,\ast}
 \| \phi_{h,\Delta t}-\psi_{h,\Delta t}\|_{0,-\frac{1}{2},\Gamma,\ast} \,,
\end{equation*}
so that we conclude:
\begin{equation*}
 \|\phi_{h,\Delta t}-\psi_{h,\Delta t}\|_{0,-\frac{1}{2},\Gamma,\ast} \lesssim||(\partial_t{f})_{h,\Delta t}-{\partial_t f}||_{0,\frac{1}{2},\Gamma} +\frac{1}{\Delta t} \|{\phi}-{\psi}_{h,\Delta t}||_{1,-\frac{1}{2},\Gamma,\ast}\,.
\end{equation*}
Using the triangle inequality, one shows that
\begin{align*}
 \|\phi-\phi_{h,\Delta t}\|_{0,-\frac{1}{2},\Gamma,\ast} &=  \|\phi-\psi_{h,\Delta t}\|_{0,-\frac{1}{2},\Gamma,\ast} +\|\phi_{h,\Delta t}-\psi_{h,\Delta t}\|_{0,-\frac{1}{2},\Gamma,\ast}\\
          &\lesssim ||(\partial_tf)_{h,\Delta t}-{\partial_tf}||_{0,\frac{1}{2},\Gamma} + \inf_{\psi_{h,\Delta t} } \{  \|{\phi}-{\psi}_{h,\Delta t}||_{0,-\frac{1}{2},\Gamma,\ast}+ \frac{1}{\Delta t}\|{\phi}-{\psi}_{h,\Delta t}||_{1,-\frac{1}{2},\Gamma,\ast}\}\\
          &\lesssim ||(\partial_tf)_{h,\Delta t}-{\partial_tf}||_{0,\frac{1}{2},\Gamma}\\
           &  \qquad \qquad + \inf_{\psi_{h,\Delta t} } \{  (1 + \frac{1}{\Delta t})\|{\phi}-{\psi}_{h,\Delta t}||_{0,-\frac{1}{2},\Gamma,\ast}+\frac{1}{\Delta t}\|\partial_t\phi-\partial_t\psi_{h,\Delta t}||_{0,-\frac{1}{2},\Gamma,\ast}\}\ .
\end{align*}
The second assertion follows from the approximation properties stated in Lemma \ref{interp}.\\

\subsection{Acoustic boundary problem}

Next, we consider the variational formulation \eqref{eq:var_time_space} of the acoustic boundary problem:\\ Find $\Phi=(\varphi,p)\in H_\sigma^1(\mathbb{R}^+,\widetilde{H}^{\frac{1}{2}}(\Gamma)) \times H_\sigma^1(\mathbb{R}^+,L^2(\Gamma))$ such that for all $\Psi=(\psi,q)\in H_\sigma^1(\mathbb{R}^+,\widetilde{H}^{\frac{1}{2}}(\Gamma)) \times H_\sigma^1(\mathbb{R}^+,L^2(\Gamma))$:
\begin{equation}\label{weaktd}
 a(\Phi,\Psi)=l(\Psi) \,.
\end{equation}

{
We obtain an a priori estimate analogous to \eqref{omegaacousticbound} in the frequency domain. If $\mathrm{Re} \ \alpha_\infty, \mathrm{Re} \ \alpha \geq 0$ and $\alpha, \frac{1}{\alpha} \in L^\infty(\Gamma)$, then for all
\begin{equation}\|p\|_{r,0,\Gamma}+ ||\varphi||_{r,\frac{1}{2},\Gamma,\ast}+  ||\varphi||_{r+1,0,\Gamma} \lesssim_\sigma \min\{|| F||_{r+1,-\frac{1}{2},\Gamma}, ||F||_{r,0,\Gamma}\} + ||G||_{r,0,\Gamma}\ .
\end{equation}
A solution in $H_\sigma^1(\mathbb{R}^+,\widetilde{H}^{\frac{1}{2}}(\Gamma)) \times H_\sigma^1(\mathbb{R}^+,L^2(\Gamma))$ therefore exists provided $G \in H_\sigma^1(\mathbb{R}^+,\widetilde{H}^{0}(\Gamma))$ and $F \in H_\sigma^2(\mathbb{R}^+,\widetilde{H}^{-\frac{1}{2}}(\Gamma))$ or $F \in H_\sigma^1(\mathbb{R}^+,\widetilde{H}^{0}(\Gamma))$.\\
}
The Galerkin discretisation of \eqref{weaktd} reads: Find $\Phi_{h,\Delta t}=(p_{h,\Delta t},\varphi_{h,\Delta t}) \in V_{h,\Delta t}^{\tilde{p},\tilde{q}} \times V_{h,\Delta t}^{p,q}$ such that for all $\Psi_{h,\Delta t}=(q_{h,\Delta t},\psi_{h,\Delta t})\in V_{h,\Delta t}^{\tilde{p},\tilde{q}} \times V_{h,\Delta t}^{p,q}$:
\begin{align}\label{eq:acoustic_Galerkin}
 a(\Phi_{h,\Delta t},\Psi_{h,\Delta t})&=\tilde{l}(\Psi_{h,\Delta t}):=\int_0^\infty\int_\Gamma  F_{h,\Delta t} \partial_t\psi_{h,\Delta t} \,ds_x\, d_\sigma t+\int_0^\infty \int_\Gamma \frac{G_{h,\Delta t} q_{h,\Delta t}}{\alpha}\, ds_x\, d_\sigma t\,.
 \end{align}

We now derive an estimate for the error of the above Galerkin approximation to \eqref{eq:var_time_space} in the norm $|||.|||_\ast$ defined by:
\begin{equation*}
|||\Phi|||_\ast=\left( ||p||_{0,0,\Gamma}^2+||\varphi||_{0,\frac{1}{2},\Gamma,\ast}^2+||\partial_t\varphi||_{0,0,\Gamma}^2\right)^{\frac{1}{2}} \quad \forall \Phi=(p,\varphi)\,.
\end{equation*}
\begin{theorem}
Assume that $\mathrm{Re} \ \alpha_\infty, \mathrm{Re} \ \alpha \geq 0$ and $\alpha, \frac{1}{\alpha} \in L^\infty(\Gamma)$.
For the solutions $\Phi=(p,\varphi)\in H_\sigma^1(\mathbb{R}^+,\widetilde{H}^{\frac{1}{2}}(\Gamma)) \times H_\sigma^1(\mathbb{R}^+,L^2(\Gamma))$ of \eqref{eq:var_time_space} and $\Phi_{h,\Delta t} =(p_{h,\Delta t},\varphi_{h,\Delta t}) \in V_{h,\Delta t}^{\tilde{p},\tilde{q}} \times V_{h,\Delta t}^{p,q}$ of \eqref{eq:acoustic_Galerkin} there holds:
\begin{align*}
&|||p-p_{h,\Delta t},\varphi-\varphi_{h,\Delta t}|||_\ast\\
&\lesssim ||F_{h,\Delta t}-F||_{0,0,\Gamma}+||G_{h,\Delta t}-G||_{0,0,\Gamma}\\
&\qquad  +  \max\left(\frac{1}{\Delta t},\frac{1}{\sqrt{h}}\right)\inf_{(q_{h,\Delta t}, \psi_{h,\Delta t})\in V_{h,\Delta t}^{\tilde{p},\tilde{q}} \times V_{h,\Delta t}^{p,q}}\left( ||p-q_{h,\Delta t}||_{1,0,\Gamma}+||\varphi-\psi_{h,\Delta t}||_{1,\frac{1}{2},\Gamma}\right)
\end{align*}
If in addition $\varphi \in H_\sigma^{s_1}(\mathbb{R}^+,{H}^{m_1}(\Gamma))$,  $p \in H_\sigma^{s_2}(\mathbb{R}^+,{H}^{m_2}(\Gamma))$,
then we have
\begin{align*}
&|||p-p_{h,\Delta t},\varphi-\varphi_{h,\Delta t}|||_\ast\\
  &\lesssim  ||F_{h,\Delta t}-F||_{0,0,\Gamma}+||G_{h,\Delta t}-G||_{0,0,\Gamma}\\
    & \qquad + \max\left(\frac{1}{\Delta t},\frac{1}{\sqrt{h}}\right)\left((h^{\alpha_1}+\Delta t^{\beta_1}) ||p||_{s_1,m_1,\Gamma}+(h^{\alpha_2}+\Delta t^{\beta_2})||\varphi||_{s_2,m_2,\Gamma}\right)\,,
\end{align*}
where
\begin{align*}
 \alpha_1 = m_1\qquad \qquad, &\quad \beta_1 = m_1+s_1-1 \,,\\
 \alpha_2 = \min\{m_2-\frac{1}{2}, m_2-\frac{3m_2}{2(m_2+s_2)}\}, & \quad \beta_2 =  m_2+s_2-\frac{3}{2}\,.
\end{align*}
\end{theorem}
\noindent \emph{Proof: } We write $\Psi=(q,\psi)$ and start with the coercivity \eqref{eq:var_coerc} applied to $\Phi_{h,\Delta t}-\Psi_{h,\Delta t} \in H^1_\sigma(\mathbb{R}^+,\widetilde{H}^{-\frac{1}{2}}(\Gamma))$ and $\Psi_{h,\Delta t}\in V_{h,\Delta t}^{\tilde{p},\tilde{q}} \times V_{h,\Delta t}^{p,q}$:
\begin{align*}
|||\Phi_{h,\Delta t}-\Psi_{h,\Delta t} |||_\ast^2&\lesssim a(\Phi_{h,\Delta t}-\Psi_{h,\Delta t},\Phi_{h,\Delta t}-\Psi_{h,\Delta t})\\
 &= a(\Phi_{h,\Delta t}-\Phi,\Phi_{h,\Delta t}-\Psi_{h,\Delta t})+a(\Phi - \Psi_{h,\Delta t},\Phi_{h,\Delta t}-\Psi_{h,\Delta t})\ .
\end{align*}
With the help of \eqref{eq:var_time_space}, \eqref{eq:bilinear_time} and \eqref{eq:rechte_time}, the first term leads to:
\begin{align*}
a(\Phi_{h,\Delta t}-\Phi,\Phi_{h,\Delta t}-\Psi_{h,\Delta t})&=\int_0^{\infty} \int_\Gamma (F_{h,\Delta t}-F)(\partial_t\varphi_{h,\Delta t}-\partial_t\psi_{h,\Delta t})\, ds_x \,d_\sigma t \\
                                                                 &\qquad + \int_0^{\infty} \int_\Gamma (\frac{G_{h,\Delta t}}{\alpha}-\frac{G}{\alpha})(p_{h,\Delta t}-q_{h,\Delta t})\, ds_x \,d_\sigma t\\
                                                                &\lesssim ||F_{h,\Delta t}-F||_{0,0,\Gamma} ||\partial_t\varphi_{h,\Delta t}-\partial_t\psi_{h,\Delta t}||_{0,0,\Gamma}\\
                                                               &\qquad \qquad + ||G_{h,\Delta t}-G||_{0,0,\Gamma} ||p_{h,\Delta t}-q_{h,\Delta t}||_{0,0,\Gamma}\\
                                                               &\leq \left(||F_{h,\Delta t}-F||_{0,0,\Gamma}+||G_{h,\Delta t}-G||_{0,0,\Gamma}\right)\\
                                 & \qquad \qquad \qquad \left(||\partial_t\varphi_{h,\Delta t}-\partial_t\psi_{h,\Delta t}||_{0,0,\Gamma}+ ||p_{h,\Delta t}-q_{h,\Delta t}||_{0,0,\Gamma}\right)\\
                                                               &\leq \left(||F_{h,\Delta t}-F||_{0,0,\Gamma}+||G_{h,\Delta t}-{G}||_{0,0,\Gamma}\right) |||\Phi_{h,\Delta t}-\Psi_{h,\Delta t}|||_\ast\,.
\end{align*}
Due to the continuity \eqref{eq:var_cont} we can estimate  the second term by
\begin{align}\label{eq:apro_1}
 |a(\Phi-\Psi_{h,\Delta t},\Phi_{h,\Delta t}-\Psi_{h,\Delta t})| \lesssim&\left( ||p-q_{h,\Delta t}||_{1,0,\Gamma}+||\varphi-\psi_{h,\Delta t}||_{1,\frac{1}{2},\Gamma,\ast}\right)\nonumber \\
                                                                                         & \,. \left( ||\varphi_{h,\Delta t}-\psi_{h,\Delta t}||_{1,\frac{1}{2},\Gamma,\ast}+||p_{h,\Delta t}-q_{h,\Delta t}||_{1,0,\Gamma}\right)\,.
\end{align}
Taking into account the inverse estimate from Section 3, we have
\begin{align}\label{eq:apro_2}
 ||\varphi_{h,\Delta t}-\psi_{h,\Delta t}||_{1,\frac{1}{2},\Gamma,\ast}&\lesssim ||\varphi_{h,\Delta t}-\psi_{h,\Delta t}||_{0,\frac{1}{2},\Gamma,\ast}+||\dot\partial_t\varphi_{h,\Delta t}-\partial_t\psi_{h,\Delta t}||_{0,\frac{1}{2},\Gamma,\ast}\nonumber\\
                                                          &\lesssim ||\varphi_{h,\Delta t}-\psi_{h,\Delta t}||_{0,\frac{1}{2},\Gamma,\ast}+(h^{-\frac{1}{2}}+\Delta t^{-\frac{1}{2}})||\partial_t\varphi_{h,\Delta t}-\partial_t\psi_{h,\Delta t}||_{0,0,\Gamma}
\end{align}
and in time
\begin{equation}\label{eq:apro_3}
 ||p_{h,\Delta t}-q_{h,\Delta t}||_{1,0,\Gamma} \leq \frac{1}{\Delta t} ||p_{h,\Delta t}-q_{h,\Delta t}||_{0,0,\Gamma} \,.
\end{equation}
Substituting \eqref{eq:apro_2} and \eqref{eq:apro_3} into \eqref{eq:apro_1} results in
\begin{align*}
 |a(\Phi-\Psi_{h,\Delta t},\Phi_{h,\Delta t}-\Psi_{h,\Delta t})| &\lesssim  (||p-q_{h,\Delta t}||_{1,0,\Gamma}+||\varphi-\psi_{h,\Delta t}||_{1,\frac{1}{2},\Gamma,\ast})\\
                                                                 & \,. \left( ||\varphi_{h,\Delta t}-\psi_{h,\Delta t}||_{0,\frac{1}{2},\Gamma,\ast}+(\frac{1}{\sqrt{h}}+\frac{1}{\sqrt{\Delta t}})||\partial_t\varphi_{h,\Delta t}-\partial_t\psi_{h,\Delta t}||_{0,0,\Gamma}\right.\\
                                                                  &\qquad \qquad \qquad \qquad \qquad\qquad \qquad \left.+\frac{1}{\Delta t}||p_{h,\Delta t}-q_{h,\Delta t}||_{0,0,\Gamma}\right)\\
                                                                 &\lesssim\, \max\left(\frac{1}{\Delta t},\frac{1}{\sqrt{h}}\right)\left( ||p-q_{h,\Delta t}||_{1,0,\Gamma}+||\varphi-\psi_{h,\Delta t}||_{1,\frac{1}{2},\Gamma}\right)\\
                                                                 &\qquad \qquad \qquad \qquad \qquad\qquad \qquad \qquad \qquad  .|||\Phi_{h,\Delta t}-\Psi_{h,\Delta t}|||_\ast\,.
\end{align*}
Altogether, we conclude
\begin{align}\label{eq:theo_4.1_inq1}
|||\Phi - \Phi_{h,\Delta t}|||_\ast &\lesssim||F_{h,\Delta t}-F||_{0,0,\Gamma}+||G_{h,\Delta t}-G||_{0,0,\Gamma}\\
                                           & \quad +\max\left(\frac{1}{\Delta t},\frac{1}{\sqrt{h}}\right)\left( ||p-q_{h,\Delta t}||_{1,0,\Gamma}+||\varphi-\psi_{h,\Delta t}||_{1,\frac{1}{2},\Gamma}\right)\nonumber\,.
\end{align}
Using the interpolation operator from Lemma \ref{interp}, we obtain the powers of $h$ and $\Delta t$ stated in the theorem.

\section{Appendix}

\noindent \emph{Proof of Theorem 4.1: }
We show that interior \eqref{eq:interior_helmholz_half} and exterior Helmholtz problems \eqref{eq:exterior_helmholz_half} with homogeneous boundary conditions $\tilde{f}=\tilde{g}=0$ admit at most one solution, $u^i = u^e = 0$. To do so we multiply the Helmholtz equation \eqref{eq:interior_helmholz_half}, \eqref{eq:exterior_helmholz_half} in $\Omega^e$ and $\Omega^i$ with $ i \bar{\omega} \bar{u}$ and integrate over $\Omega^e \cup \Omega^i$. We obtain
\begin{equation*}
 \int_{\Omega^e \cup \Omega^i} \Delta u\cdot i\bar{\omega} \bar{u}+\omega^2 u\cdot i\bar{\omega} \bar{u} \, dx=0\,.
\end{equation*}
Applying Green's first theorem to $u^e$ and $u^i$, we obtain
\begin{equation*}
\int_{\Omega^e \cup \Omega^i} -i\bar{\omega}|\triangledown u|^2 +i \omega |\omega|^2 |u|^2  \, dx -\int_\Gamma i\bar{\omega}(\frac{\partial u^e}{\partial n} \bar{u}^e -\frac{\partial u^i}{\partial n} \bar{u}^i)\, ds_x-\int_{\Gamma_\infty\cup \Gamma_\infty'} i\bar{\omega}\frac{\partial u}{\partial n} \bar{u}\, ds_x=0\,.
\end{equation*}
Here, we have neglected a contribution from a large half--sphere which tends to zero as the radius of the half--sphere goes to infinity.\\
We take the real part of this equality and use the boundary conditions:
\begin{align*}
 &2\mathrm{Im}\ \omega \int_{\Omega^e \cup \Omega^i} |\triangledown u|^2 + |\omega|^2 |u|^2 \, dx = \mathrm{Re}(\int_\Gamma -i\bar{\omega}(\frac{\partial u^e}{\partial n} \bar{u}^e -\frac{\partial u^i}{\partial n} \bar{u}^i)\,ds_x-\int_{\Gamma_\infty\cup \Gamma_\infty'} i\bar{\omega} \frac{\partial u}{\partial n} \bar{u}\,ds_x)\\
&= \mathrm{Re}(\int_\Gamma -i\bar{\omega}(-\alpha i \omega u^e \bar{u}^e -\alpha i \omega u^i \bar{u}^i)\,ds_x-i\bar{\omega} \int_{\Gamma_\infty\cup \Gamma_\infty'} -\alpha_\infty i \omega u \bar{u}\,ds_x)\\
&=  -(\mathrm{Re}\ \alpha) \int_\Gamma  |\omega|^2 |u^e|^2 + |\omega|^2 |u^i|^2 \,ds_x -(\mathrm{Re}\ \alpha_\infty) \int_{\Gamma_\infty\cup \Gamma_\infty'} |\omega|^2 |u|^2 \,ds_x\,.
\end{align*}
Since $\mathrm{Im}\ \omega > 0$, the conditions $\mathrm{Re}\ \alpha\geq 0$ and $\mathrm{Re}\ \alpha_\infty\geq 0$ ensure that $u=0$ in $H^1(\Omega^e\cup\Omega^i)$. The uniqueness of the solution follows.\\

\noindent \emph{Proof of Theorem \ref{3.4}: }
First we prove \eqref{eq:cont_V}.\\
Let be $p$ in $\widetilde{H}^{-\frac{1}{2}}(\Gamma)$ and  let $v=S_\omega p$. Then we saw that $v$ verifies:
\begin{equation*}
\begin{cases}
(\Delta +\omega^2) v(x)=0 \quad\mbox{ in }\Omega^i \cup \Omega^e\\
\frac{\partial v^i}{\partial n} -\frac{\partial v^e}{\partial n} = p \quad\mbox{ in }\Gamma \\
v^i -v^e = 0 \quad\mbox{ in }\Gamma \,.
\end{cases}
\end{equation*}
Applying Green's Theorem in $\Omega^i$ resp.~$\Omega^e$ we obtain
\begin{equation}\label{eq:innen_green}
-\int_{\Gamma_\infty'} i\bar{\omega} \frac{\partial v^i}{\partial x_3} \bar{v^i}\,ds_x+ \int_\Gamma i\bar{\omega} \frac{\partial v^i}{\partial n} \bar{v^i}\,ds_x=\int_{\Omega^i} -i\omega |\omega|^2 v^i \bar{v^i}\,dx+\int_{\Omega^i} i\bar{\omega} \triangledown v^i \triangledown \bar{v^i}\,dx
\end{equation}
and
\begin{equation}\label{eq:aussen_green}
-\int_{\Gamma_\infty\cup\Gamma_\infty'} i\bar{\omega} \frac{\partial v}{\partial x_3} \bar{v}\,ds_x-\int_\Gamma i\bar{\omega} \frac{\partial v^e}{\partial n} \bar{v^e} \,ds_x=\int_{\Omega^e} -i\omega |\omega|^2 v^e \bar{v^e} \,dx +\int_{\Omega^e} i\bar{\omega} \triangledown v^e \triangledown \bar{v^e} \,dx\,.
\end{equation}
Adding the two equations \eqref{eq:innen_green} and \eqref{eq:aussen_green} we get
\begin{equation*}
-\int_{\Gamma_\infty\cup \Gamma_\infty'} i\bar{\omega} \frac{\partial v}{\partial x_3} \bar{v} \,ds_x+\int_\Gamma i\bar{\omega} \left( \frac{\partial v^i}{\partial n}- \frac{\partial v^e}{\partial n}\right) \bar{v^e}\,ds_x=\int_{\Omega^e \cup \Omega^i} -i\omega |\omega|^2 |v|^2 \,dx +\int_{\Omega^e \cup \Omega^i} i\bar{\omega} |\triangledown v|^2 \,dx.
\end{equation*}
Using the boundary conditions on $\Gamma$ and $\partial\mathbb{R}^3_+$ we obtain
\begin{equation*}
-\int_{\Gamma_\infty\cup \Gamma_\infty'} \alpha_\infty |\omega|^2 |v|^2 \,ds_x +\int_\Gamma i\bar{\omega} p \bar{v^e}\,ds_x=\int_{\Omega^e \cup \Omega^i} -i\omega |\omega|^2 |v|^2\,dx+\int_{\Omega^e \cup \Omega^i} i\bar{\omega} |\triangledown v|^2 \,dx\,.
\end{equation*}
We take the real part of this equation:
\begin{equation*}
 -\int_{\Gamma_\infty\cup \Gamma_\infty'} \mathrm{Re}(\alpha_\infty)|\omega|^2 |v|^2 \,ds_x+ \mathrm{Re}\left(\int_\Gamma i\bar{\omega} p \bar{v^e} \,ds_x \right)=\sigma\left(||v^e||_{1,\omega,\Omega^e}^2+ ||v^i||_{1,\omega,\Omega^i}^2\right)\,.
\end{equation*}
It follows from $\mathrm{Re}(\alpha_\infty)\geq 0$ and from the trace theorem (Lemma 1.4) that
\begin{equation*}
 \mathrm{Re}\left(\int_\Gamma i\bar{\omega} p \bar{v^e}\,ds_x \right) \geq \sigma \left(||v^e||_{\frac{1}{2},\omega,\Gamma}^2+ ||v^i||_{\frac{1}{2},\omega,\Gamma}^2\right)\,.
\end{equation*}
Therefore $|\omega| ||p||_{-\frac{1}{2},\omega,\Gamma,\ast} ||v^e||_{\frac{1}{2},\omega,\Gamma} \geq 2 \sigma||v^e||_{\frac{1}{2},\omega,\Gamma}^2
$, or  $|\omega| ||p||_{-\frac{1}{2},\omega,\Gamma,\ast} \geq 2 \sigma||v^e||_{\frac{1}{2},\omega,\Gamma}$.
As $v|_{\Gamma}=V_\omega p$ we obtain
\begin{equation*}
 ||V_\omega p||_{\frac{1}{2},\omega,\Gamma} \leq \frac{|\omega|}{2\sigma} ||p||_{-\frac{1}{2},\omega,\Gamma,\ast}\,.
\end{equation*}
We now consider the estimate \eqref{eq:cont_W}. Let $\varphi$ in $\widetilde{H}^{\frac{1}{2}}(\Gamma)$ and  let $v=-D_\omega \varphi$. Then we have seen that $v$ verifies:
\begin{equation*}
\begin{cases}
(\Delta +\omega^2) v(x)=0 \quad\mbox{ in }\Omega^i \cup \Omega^e\\
\frac{\partial v^i}{\partial n} -\frac{\partial v^e}{\partial n} = 0 \quad\mbox{ in }\Gamma \\
v^i -v^e = \varphi \quad\mbox{ in }\Gamma\,.
\end{cases}
\end{equation*}
Moreover, $\frac{\partial v^e}{\partial n}|_\Gamma=-W_\omega \varphi$.\\
Adding \eqref{eq:innen_green} and \eqref{eq:aussen_green}, we obtain
\begin{equation*}
-\int_{\Gamma_\infty\cup \Gamma_\infty'} i\bar{\omega} \frac{\partial v}{\partial x_3} \bar{v}+\int_\Gamma i\bar{\omega}  \frac{\partial v^e}{\partial n} \left(\bar{v^i}-\bar{v^e}\right)=\int_{\Omega^e \cup \Omega^i} -i\omega |\omega|^2 |v|^2+\int_{\Omega^e \cup \Omega^i} i\bar{\omega} |\triangledown v|^2 \,.
\end{equation*}
Using the boundary condition on $\Gamma$ and $\Gamma_\infty\cup \Gamma_\infty'$ leads to the following equality:
\begin{equation*}
-\int_{\Gamma_\infty\cup\Gamma_\infty'} \alpha_\infty |\omega|^2 |v|^2 \,ds_x+\int_\Gamma i\bar{\omega} \frac{\partial v^e}{\partial n}\bar{\varphi}\,ds_x =\int_{\Omega^e \cup \Omega^i} -i\omega |\omega|^2 |v|^2\,dx+\int_{\Omega^e \cup \Omega^i} i\bar{\omega} |\triangledown v|^2 \,dx\,.
\end{equation*}
Its real part is given by
\begin{equation*}
 -(\mathrm{Re}\ \alpha_\infty) \int_{\Gamma_\infty\cup\Gamma_\infty'} |\omega|^2 |v|^2 \,ds_x + \mathrm{Re}\left(\int_\Gamma i\bar{\omega}  \frac{\partial v^e}{\partial n} \bar{\varphi} \,ds_x\right)=\sigma\left(||v^e||_{1,\omega,\Omega^e}^2+ ||v^i||_{1,\omega,\Omega^i}^2\right)\,.
\end{equation*}
As $\mathrm{Re}\ \alpha_\infty\geq 0$ and using Cauchy-Schwarz, we conclude
\begin{equation}\label{eq:ve_leq_ven}
 ||v^e||_{1,\omega,\Omega^e}^2 \leq \frac{1}{\mathrm{Im}(\omega)}|\omega| ||\varphi||_{\frac{1}{2},\omega,\Gamma,\ast} ||\frac{\partial v^e}{\partial n}||_{-\frac{1}{2},\omega,\Gamma}\,.
\end{equation}
It remains to estimate $\frac{\partial v^e}{\partial n}$. From Green's theorem in $\Omega^e$ we see that
\begin{equation*}
-\int_{\Gamma_\infty}  \frac{\partial v^e}{\partial x_3} \bar{\varphi} \,ds_x -\int_\Gamma  \frac{\partial v^e}{\partial n} \bar{\varphi}\,ds_x=\int_{\Omega^e} -\omega^2 v^e \bar{\psi}\,dx+\int_{\Omega^e}  \triangledown v^e \triangledown \bar{\psi}\,dx\,,
\end{equation*}
where $\quad \psi \in H^1(\Omega^e)$ is an extension of $\phi \in \widetilde{H}^{\frac{1}{2}}(\Gamma)$ with $||\psi||_{1,\omega,\Omega^e} \lesssim_\sigma ||\phi||_{\frac{1}{2},\omega,\Gamma,\ast}$.\\
From the trace theorem
\begin{equation*}
 |\int_\Gamma  \frac{\partial v^e}{\partial n} \bar{\varphi}\,ds_x| \leq ||v^e||_{1,\omega,\Omega^e} ||\psi||_{1,\omega,\Omega^e}\,.
\end{equation*}
we conclude
\begin{equation*}
 ||\frac{\partial v^e}{\partial n}||_{-\frac{1}{2},\omega,\Gamma}=\sup_{\{\phi \in \widetilde{H}^{\frac{1}{2}}(\Gamma) / ||\phi||_{\frac{1}{2},\omega,\Gamma,\ast}=1\}} |\int_\Gamma  \frac{\partial v^e}{\partial n} \bar{\phi}\,ds_x| \leq C ||v^e||_{1,\omega,\Omega^e}\,.
\end{equation*}
From \eqref{eq:ve_leq_ven} it follows that
\begin{equation*}
 ||W_\omega \varphi||_{-\frac{1}{2},\omega,\Gamma}=||\frac{\partial v^e}{\partial n}||_{-\frac{1}{2},\omega,\Gamma} \leq C ||v^e||_{1,\omega,\Omega^e} \leq C|\omega| ||\varphi||_{\frac{1}{2},\omega,\Gamma,\ast}\,.
\end{equation*}
Using similar reasoning, we obtain the estimates \eqref{eq:cont_K'} and \eqref{eq:cont_K}.\\

\noindent \emph{Proof of Theorem \ref{3.5}: }
Recall that with $\tilde{U}=(\varphi,p)$ and $\tilde{V}=(\psi,q)$
\begin{align*}\label{eq:bilinear_freq}
a_\omega(\tilde{U},\tilde{V})=|\omega|^2 \int_\Gamma \alpha \varphi \bar{\psi} ds_x + \int_\Gamma \frac{1}{\alpha} p \bar{q} ds_x +  i \bar{\omega}\int_\Gamma K'_\omega p \bar{\psi} ds_x\\
-i \bar{\omega}\int_\Gamma W_\omega \varphi \bar{\psi} ds_x - i \omega \int_\Gamma V_\omega p \bar{q} ds_x + i \omega \int_\Gamma K_\omega \varphi \bar{q} ds_x\,.
\end{align*}
We estimate the various terms of the bilinear form $a_\omega$ using Theorem \ref{3.4}
\begin{equation*}
 |\omega|^2 \left|\int_\Gamma \alpha \varphi \bar{\psi} ds_x\right|\leq C |\omega|^2 ||\varphi||_{0,\omega,\Gamma} ||\psi||_{0,\omega,\Gamma} \leq C |\omega|^2 ||\varphi||_{\frac{1}{2},\omega,\Gamma,\ast} ||\psi||_{\frac{1}{2},\omega,\Gamma,\ast}\,,
\end{equation*}
\begin{equation*}
\left|\int_\Gamma \frac{1}{\alpha} p \bar{q} ds_x \right| \leq C  ||p||_{0,\omega,\Gamma} ||q||_{0,\omega,\Gamma}\leq \frac{C}{\sigma^2} |\omega|^2 ||p||_{0,\omega,\Gamma} ||q||_{0,\omega,\Gamma}\,,
\end{equation*}
\begin{align*}
 \left| i \bar{\omega}\int_\Gamma K'_\omega p \bar{\psi} ds_x\right| &\leq \frac{|\omega|}{\sigma} ||K'_\omega p||_{-\frac{1}{2},\omega,\Gamma} ||\psi||_{\frac{1}{2},\omega,\Gamma,\ast}\\
                                                                        &\leq C\frac{|\omega|^2}{\sigma} ||p||_{-\frac{1}{2},\omega,\Gamma,\ast} ||\psi||_{\frac{1}{2},\omega,\Gamma,\ast}\\
                                                                        &\leq C\frac{|\omega|^2}{\sigma} ||p||_{0,\omega,\Gamma} ||\psi||_{\frac{1}{2},\omega,\Gamma,\ast}\,,
\end{align*}
\begin{equation*}
 \left| i \bar{\omega}\int_\Gamma W_\omega \varphi \bar{\psi} ds_x\right| \leq |\omega| ||W_\omega \varphi||_{-\frac{1}{2},\omega,\Gamma} ||\psi||_{\frac{1}{2},\omega,\Gamma,\ast} \leq C|\omega|^2 ||\varphi||_{\frac{1}{2},\omega,\Gamma,\ast} ||\psi||_{\frac{1}{2},\omega,\Gamma,\ast}\,,
\end{equation*}
\begin{align*}
 \left| i \omega \int_\Gamma V_\omega p \bar{q} ds_x\right| &\leq |\omega| ||V_\omega p||_{\frac{1}{2},\omega,\Gamma} ||q||_{-\frac{1}{2},\omega,\Gamma,\ast}\\
                                                                & \leq C|\omega|^2 ||p||_{-\frac{1}{2},\omega,\Gamma,\ast} ||q||_{-\frac{1}{2},\omega,\Gamma,\ast}\\
                                                                &\leq C|\omega|^2 ||p||_{0,\omega,\Gamma} ||q||_{0,\omega,\Gamma}
\end{align*}
\begin{align*}
 \left| i \omega \int_\Gamma K_\omega \varphi \bar{q} ds_x\right| &\leq |\omega| ||K_\omega \varphi||_{\frac{1}{2},\omega,\Gamma} ||q||_{-\frac{1}{2},\omega,\Gamma,\ast}\\
                                                                        &\leq C|\omega|^2 ||\varphi||_{\frac{1}{2},\omega,\Gamma,\ast} ||q||_{-\frac{1}{2},\omega,\Gamma,\ast}\\
                                                                        &\leq C|\omega|^2 ||\varphi||_{\frac{1}{2},\omega,\Gamma,\ast} ||q||_{0,\omega,\Gamma}
\end{align*}
Adding the 6 inequalities we get
\begin{equation}\label{eq:a_freq_cont}
a_\omega(\tilde{U},\tilde{V}) \lesssim_\sigma \left({|\omega|||p||}_{0,\omega,\Gamma} +|\omega|||\varphi||_{\frac{1}{2},\omega,\Gamma,\ast}\right)\left(|\omega|||q||_{0,\omega,\Gamma} +|\omega|||\psi||_{\frac{1}{2},\omega,\Gamma,\ast}\right)\,.
\end{equation}


\begin{thebibliography}{99}
\bibitem{bamberger86}A.~Bamberger and T.~Ha Duong, \emph{Formulation variationnelle espace-temps pour le calcul par potentiel retardé de la diffraction d'une onde acoustique}, Math.~Meth.~Appl.~Sci.~\textbf{8} (1986), 405--435 and 598--608.

\bibitem{terrasse93}I.~Terrasse, \emph{R\'{e}solution math\'{e}matique et num\'{e}rique des \'{e}quations de Maxwell instationnaires par une m\'{e}thode de potentiels retard\'{e}s}, PhD thesis, \'{E}cole Polytechnique, Palaiseau, 1993.

\bibitem{hld03} T.~Ha-Duong, B.~Ludwig and I.~Terrasse, {\em A Galerkin BEM for transient acoustic scattering by an absorbing obstacle}, Internat.~J.~Numer.~Methods Engrg. \textbf{57} (2003), 1845--1882.

\bibitem{mic} A.~E.~Yilmaz, J.-M.~Jin, E.~Michielssen, \emph{Time domain adaptive integral method for surface integral equations}, IEEE Trans.~Antennas Propagation \textbf{52} (2004) 2692--2708.

\bibitem{dd1} P.~J.~Davies, D.~B.~Duncan, \emph{Convolution-in-time approximations of time domain boundary integral equations}, SIAM J.~Sci.~Comput.~\textbf{35} (2013), B43--B61.

\bibitem{dd2} P.~J.~Davies, D.~B.~Duncan, \emph{Convolution spline approximations for time domain boundary integral equations}, J.~Integral Equations Applications, to appear (2014).

\bibitem{sayas} F.~J.~Sayas, \emph{Retarded Potentials and Time Domain Boundary Integral Equations: a road-map}, lecture notes, 2013.

\bibitem{Ochmann01} {M.~Ochmann}, {\em Closed form solutions for the acoustical impulse response over a masslike or an absorbing plane}, {J. Acoust. Soc. Am.} {129 (6)}, {2011}.

\bibitem{Ha-Duong03a} T.~Ha Duong, \emph{On retarded potential boundary integral equations and their discretisations}, in: Topics in computational wave propagation, pp.~301–-336, Lect.~Notes Comput.~Sci.~Eng., \textbf{31}, Springer, Berlin, 2003.

\bibitem{costabel04}
{M.~Costabel}, {\em Time-dependent problems with the boundary integral equation method}. {In Encyclopedia of Computational Mechanics}, {E. ~Stein, ~R.~de
~Borst, ~and ~J. ~R. Hughes}, {Eds. John Wiley \& Sons}, {Chichester}, {2004},{pp. 703-721}.

\bibitem{Gimperlein} H.~Gimperlein, Z.~Nezhi, E.~P.~Stephan, \emph{A residual a posteriori error estimate for the time--domain boundary element method}, in preparation.

\bibitem{Banz} L.~Banz, H.~Gimperlein, Z.~Nezhi, E.~P.~Stephan, \emph{Time domain BEM for sound radiation of tyres}, in preparation.

\bibitem{tri}H.~Triebel, \emph{Theory of Function Spaces I/II}, Birkh\"auser, Basel, 1983/1992.

\bibitem{screen} {E.~P.~Stephan}, {\em Boundary integral equations for screen problems in $\mathbb{R}^3$}, {Integral Equations Operator Theory} {10} (1987), 236 -- 257.

\bibitem{Heuer} N.~Heuer, \emph{Additive Schwarz method for the $p$--version of the boundary element method for the single layer potential operator on a plane screen}, Numer.~Math.~\textbf{88} (2001), 485--511.

\bibitem{Glaefke} {M.~Glaefke}, {\em Adaptive Methods for Time Domain Boundary Integral Equations}. {PhD
thesis, Brunel University, 2012}.

\bibitem{Babuska_aziz}I.~Babuska and A.~K.~Aziz, \emph{Survey lectures on the mathematical foundations of the finite element method}. In: The Mathematical Foundations of the Finite Element Method with Applications to Partial Differential Equations, Academic Press, New York, 1972, 3-359.

\bibitem{Ochmann02} {M.~Ochmann}, {\em The complex equivalent source method for sound propagation over an impedance plane}, {J. Acoust. Soc. Am.} {116 (6)}, {2004}.

\bibitem{becache94} {E.~B\'{e}cache}, {\em Equations int\'{e}grales pour l'\'{e}quation des ondes}. {Cours de l'\'{e}cole des ondes}, {INRIA}, {1994}.
{Available for download at http://www-rocq.inria.fr/~becache/cours\_eqinteg.ps.gz}.

\bibitem{nedelec} J.~C.~Nedelec, \emph{Curved finite element methods for the solution of singular integral equations on surfaces in $\mathbb{R}^3$},  Comput.~Methods Appl.~Mech.~Engrg.~\textbf{8} (1976), 61–-80.



\end{thebibliography}
\end{document}